\documentclass[10 pt, twocolumn, twoside]{IEEEtran}  	

\IEEEoverridecommandlockouts                              

\usepackage{times}
\usepackage{url}
\usepackage{graphicx}
\usepackage{subfig}
\usepackage{epsfig} 
\usepackage{times} 
\usepackage{amsmath} 
\usepackage{amssymb}  
\usepackage{euscript}
\usepackage[usenames,dvipsnames]{xcolor}
\usepackage{algorithmic}
\usepackage{algorithm}
\usepackage{dsfont}
\usepackage{bbm}
\usepackage{enumerate}
\usepackage[subnum]{cases}
\graphicspath{{./images/}}
\usepackage{caption}
\usepackage{tikz}
\usepackage{pgfplots}

\DeclareSymbolFont{extraup}{U}{zavm}{m}{n}
\DeclareMathSymbol{\vardiamond}{\mathalpha}{extraup}{87}

\usepackage[normalem]{ulem}

\newtheorem{theorem}{Theorem}

\newtheorem{definition}[theorem]{Definition}
\newtheorem{proposition}[theorem]{Proposition}
\newtheorem{assumption}[theorem]{Assumption}

\newtheorem{lemma}[theorem]{Lemma}
\newtheorem{remark}[theorem]{Remark}
\newtheorem{problem}[theorem]{Problem}
\newtheorem{example}[theorem]{Example}

\def\R{{\mathbb{R}}}

\def\V{{\mathbb{V}}}

\def\Qcrit{Q_{\rm crit}}
\def\Qcritinv{\Qcrit^{-1}}
\def\1{\mathbbm{1}}

\newcommand{\until}[1]{\{1,\dots, #1\}}

\def\secure_set{\V_\alpha}

\usepackage{color}

\definecolor{forestgreen}{RGB}{32,178,170}
\definecolor{golden}{RGB}{255,215,0}
\newcommand{\mycbox}[1]{\tikz{\path[draw=#1,fill=#1] (0,0) rectangle (2mm,2mm);}}

\newcommand{\tq}[1]{\textquotedblleft \textnormal{#1}\textquotedblright}

\newcommand{\vect}[1]{\mathbbold{#1}}
\newcommand{\vones}[1][]{\vect{1}_{#1}}
\newcommand{\vzeros}[1][]{\vect{0}_{#1}}
\DeclareSymbolFont{bbold}{U}{bbold}{m}{n}
\DeclareSymbolFontAlphabet{\mathbbold}{bbold}
\newcommand\oprocendsymbol{\hbox{$\square$}}
\newcommand\oprocend{\relax\ifmmode\else\unskip\hfill\fi\oprocendsymbol}

\title{\LARGE\bf
Voltage stress minimization by optimal reactive power control
}

\author{Marco~Todescato, John~W.~Simpson-Porco, Florian~D\"orfler, Ruggero~Carli and Francesco~Bullo
  \thanks{This work was supported by the Ing. Aldo Gini Foundation, Padova and by the ETH start-up funds.
   M.~Todescato and R.~Carli are with the
   Department of Information Engineering, University of Padova 
   {\tt\footnotesize \{todescat|carlirug\}
   @dei.unipd.it}.
   J.~W.~Simpson-Porco is with the Department of Electrical and Computer Engineering, University of Waterloo {\tt\footnotesize jwsimpson@uwaterloo.ca}. F.~Bullo is with the Department of Mechanical
   Engineering and the Center for Control, Dynamical System and Computation,
   University of California at Santa Barbara
   {\tt\footnotesize bullo@engineering.ucsb.edu}.
   F.~D\"orfler is with the Automatic Control Laboratory, Swiss Federal Institute (ETH) Zurich
   {\tt\footnotesize dorfler@ethz.ch}.}}

\begin{document}

\maketitle
\thispagestyle{empty}
\pagestyle{empty}

\begin{abstract}
A standard operational requirement in power systems is that the voltage magnitudes  lie within prespecified bounds. 
Conventional engineering wisdom suggests that such a tightly-regulated profile, imposed for system design purposes and good operation of the network, should also guarantee a secure system, operating far from static bifurcation instabilities such as voltage collapse. 
%
In general however, these two objectives are distinct and must be separately enforced. 
We formulate an optimization problem which maximizes the distance to voltage collapse through injections of reactive power, subject to power flow and operational voltage constraints.
By exploiting a linear approximation of the power flow equations we arrive at a convex reformulation which can be efficiently solved for the optimal injections. 
We also address the planning problem of allocating the resources by recasting our problem in a sparsity-promoting framework that allows us to choose a desired trade-off between 
optimality of injections and the number of required actuators.
%
%
Finally, we present a distributed algorithm to solve the optimization problem, showing that it can be implemented on-line as a feedback controller.
We illustrate the performance of our results with the IEEE30 bus network.
\end{abstract}

\begin{keywords}
\small power networks, voltage support, reactive power compensation, resource allocation, distributed control. 
\end{keywords}

\section{Introduction}\label{sec:intro}
%
Traditionally, the main purpose of voltage support is to maintain voltage magnitudes tightly within predetermined security constraints (e.g., within $5\%$ of some nominal level). Conventional wisdom suggests that such a tightly regulated voltage profile, imposed for system design reasons and good operation of the power network, should also guarantee a secure system, operating far from static bifurcation instabilities such as voltage collapse.
%
%
Techniques for voltage support include shunt and static VAR compensation  \cite{ASY-NM:04}, 
series compensation \cite{MM-MF:00}, off-nominal transformer tap ratios \cite{FAV-AS-JD:07}, synchronous condensers \cite{GA:08}, and inverters operating away from unity power factor \cite{NL-GQ-MD:14}. 
See \cite{JD-LM-LR-RD:05} for a survey on the topic. 

A distinct voltage control problem, which represents a key direction in power system stability analysis, has been the development of indices quantifying a power network's proximity to voltage collapse. A broad overview of this large subfield can be found in 
\cite{TVC:00,CC:02-short,TVC-CV:98}. 
The most reliable existing approaches are largely based on numerical methods and lack detailed theoretical support. They often require either continuation power flow \cite{IAH-RJD:01} to identify the insolvability boundary, or repeated computation of loading margins over various directions in parameter-space \cite{ID-LL:92b}.

As stressed, voltage support and distance to collapse are often analyzed separately in power systems although they are intrinsically related through the well known principle of reactive power injection. Combining the two problems represents the first contribution of the paper which is threefold.
Indeed the ultimate goal in voltage support problems is the \textit{security task} to confine the voltage magnitudes within predetermined bounds, as suggested by conventional engineering wisdom. Here, we follow an alternative approach: we define a particular measure for the \textit{network stress}, i.e., the stress experienced by the network induced by the load profile.
In particular, we begin our analysis from the recent article \cite{JWSP-FD-FB:14c} where a sufficient and tight condition was presented for solvability of \emph{decoupled reactive power flow}. This condition, rigorously proved only for the reactive decoupled case, quantifies the proximity to voltage collapse by determining a nodal measure of network stress.
Based on this condition we pursue a novel system-level formulation of optimal voltage support encoded as an optimization problem with \emph{stress-minimization}, i.e., maximization of the distance to voltage collapse, as objective and subject to voltage security constraints. This approach allows us to match a \emph{local} security requirement as well as a \emph{system-level} stress-minimization objective encoding the distance to collapse.
By exploiting an opportune linearized reformulation, our optimization formulation becomes convex and can be efficiently solved for the optimal injections.
As second contribution, we also address the planning problem of allocating the available resources by regularizing our optimization problem with a convex proxy of the cardinality function. This sparsity-promoting formulation allows us to choose a desired trade-off between performance and a cost-effective solution.
Finally, we present a distributed algorithm for the stress
minimization problem which is amenable to real-time implementation as a
distributed feedback controller.

Compared to other approaches to voltage support problems our results do not rely on the assumption of a radial (i.e., acyclic) power grid topology \cite{NL-GQ-MD:14}. 
This makes our approach appealing for power transmission networks. Different from the reactive power compensation literature 
\cite{SB-RC-GC-SZ:15} 
and from the voltage support literature \cite{NL-GQ-MD:14}, we seek stress minimization rather than optimal power flow (minimizing, e.g., losses) or voltage security tasks. Moreover, our formulation can nicely incorporate controller placement tasks.

The remainder of this paper is organized as follows. In Section \ref{sec:preliminaries} we introduce the required power system model. In Section \ref{sec:problem} we present the first two main contributions of the paper: (i) we review the typical objectives for voltage regulation problems, propose a novel measure for the network stress, and formulate our optimization problem. We then present and solve the convex reformulation of the problem. (ii) We analyze a sparsity-promoting cost to address the planning problem. In Section \ref{sec:distributed_algorithm} we present the third  contribution consisting in a distributed strategy to perform real-time stress minimization. Finally, Section \ref{sec:conclusions} offers conclusions and future directions.

\section{Preliminaries}\label{sec:preliminaries}

\subsection{Power Network, Generator and Load Models}\label{subsec:model}
A high voltage power network can be modeled as a connected, undirected and complex-weighted graph $\mathcal{G}(\mathcal{V},\mathcal{E})$ where $\mathcal{V}= \until{n}$ represents the set of nodes (or \emph{buses}), $\mathcal{E}$ ($|\mathcal{E}| = m$) is the set of edges (or \emph{branches}) connecting the nodes, that is the set of unordered pairs $(h,k),\ h,k\in\mathcal{V}$, such that $h$ and $k$ are connected to each other. 
Under synchronous steady-state operating conditions, all the electric quantities are sinusoidal signals at the same frequency. At every bus $h \in\mathcal{V}$ we have the following phasor 
quantities:
\begin{itemize}
\item nodal voltage: $u_h = V_h\exp(\mathrm{j}\theta_h)\in\mathbb{C}$;
\item current injection: $i_h = I_h\exp(\mathrm{j}\psi_h)\in\mathbb{C}$;
\item power injection: $s_h = p_h + \mathrm{j}q_h = v_h\overline{i_h}\in\mathbb{C}$;.
\end{itemize}
where $V_h, \theta_h, I_h, \psi_h, p_h, q_h \in\mathbb{R}$ and $\overline{(\cdot)}$ denotes the complex conjugate operator. 
By collecting all quantities
into vectors $u,\ i,\ s\in\mathbb{C}^n$, Kirchhoff's and Ohm's laws lead to
\begin{equation}\label{eq:kirchhoff_current}
i = \mathrm{j}Bu\,,
\end{equation}
where $\mathrm{j}$ denotes the imaginary unit.  
The symmetric and sparse susceptance matrix $B \in \R^{n\times n}$ encodes the topology of the underlying electric network weighted by the line susceptances. Following standard assumptions we neglect line losses in high-voltage transmission networks \cite{ARB-VV:06,JM-JWB-JRB:08}. For recent studies on lossy distribution networks, see~\cite{NL-GQ-MD:14,SB-RC-GC-SZ:15}. 
From \eqref{eq:kirchhoff_current} we can write the \emph{Power Flow Equations} (PFEs) as
\begin{equation}\label{eq:comples_power_vector}
s = [u]\overline{i} = [u]\overline{(\mathrm{j}Bu)}\,,
\end{equation}
where $[x]$ denotes the diagonal $n\times n$ matrix with diagonal entries $x_i$.
Expanding \eqref{eq:comples_power_vector}, for each $h\in\mathcal{V}$ the real and imaginary parts must satisfy
\begin{subequations}
\begin{equation}\label{eq:active_pfes}
p_h(u) = \sum_{k=1}^n \nolimits B_{hk} V_hV_k\sin(\theta_h-\theta_k)\,,
\end{equation}
\begin{equation}\label{eq:reactive_pfes}
q_h(u) = -\sum_{k=1}^n \nolimits B_{hk} V_hV_k\cos(\theta_h-\theta_k)\,. 
\end{equation}
\end{subequations}
The PFEs \eqref{eq:active_pfes}--\eqref{eq:reactive_pfes} relate the voltage variables $(\theta,V)$ to the power variables $(p,q)$, while the behavior of each bus is specified by the particular model assumed to describe it. In this paper, we partition the set of buses $\mathcal{V}$ into two subsets, namely $\mathcal{V}_L$ ($|\mathcal{V}_L| = n_\ell$) which identifies \emph{power-regulated} or \emph{load} buses, and $\mathcal{V}_G$ ($|\mathcal{V}_G| = n_g$) which identifies \emph{voltage-regulated} buses.\footnote{We use the subscript $G$ for voltage-regulated buses because in transmission grids these are typically generator buses.} In particular, we assume the following:
\begin{itemize}
\item \emph{Voltage-regulated bus model}: voltage-regulated buses are modeled as standard $PV$ buses \cite{ARB-VV:06}.
This model is widely used, e.g., for generators in transmission grids and micro-generators used for reactive power control in (micro) distribution grids \cite{RHL:02}.
\item \emph{Load bus model}: loads are modeled as $PQ$ buses \cite{JM-JWB-JRB:08,MKP:92}.
In our setup, this model refers also to sources interfaced with power electronics and voltage support equipment such as synchronous condensers. While our results extend to ZIP load models \cite{ARB-VV:06}, for simplicity of presentation we restrict ourselves to constant power loads $q_h(u) = Q_h$; constant impedance loads can be incorporated into the $B$ matrix as diagonal elements.
\end{itemize}

After relabeling the buses to place loads before generators, the matrix $B$ can be partitioned in the block-matrix form
\begin{equation}\label{eq:Bmatrix}
B = \begin{pmatrix}
B_{LL} & B_{LG}\\ B_{GL} & B_{GG}
\end{pmatrix}\,.
\end{equation}

\begin{assumption}[Properties of $B_{LL}$]\label{ass:Mmatrix}~ 
\begin{enumerate}[(i)]
\item $B_{LL}$ is a Metzler matrix whose eigenvalues are characterized by a negative real part;\footnote{In other words, $B_{LL}$ is an $M$-matrix.}
\item the graph associated to the $B_{LL}$ matrix (i.e., the graph induced by the load buses $\mathcal{V}_L$) is connected.
\end{enumerate}
\end{assumption}

\smallskip

Assumption \ref{ass:Mmatrix} (i) is typically verified in practice \cite{JT-DS-MIS:86}, and always satisfied in the absence of phase-shifting transformers, line-charging and shunt capacitors. Regarding the shunt capacitors, they are allowed to be different from zero, however Assumption \ref{ass:Mmatrix} (i) limits their sizes. Assumption \ref{ass:Mmatrix} (ii) can be made without loss of generality, since connected components of the induced graph will be electrically isolated from one another by voltage-regulated generator buses.

\begin{figure*}[t]
\centering
\subfloat[Example system.\label{subfig:2node_case}]{\includegraphics[scale=0.9]{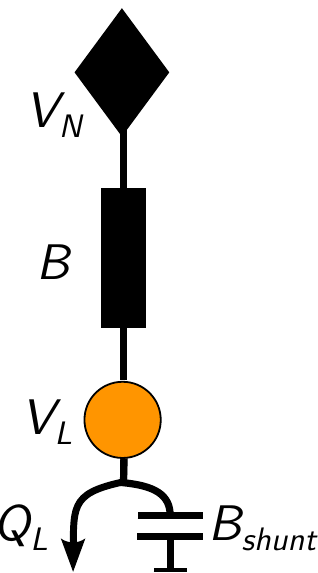}}
\subfloat[Absence of shunt.\label{subfig:no_shunts}]{\includegraphics[width=0.4\textwidth]{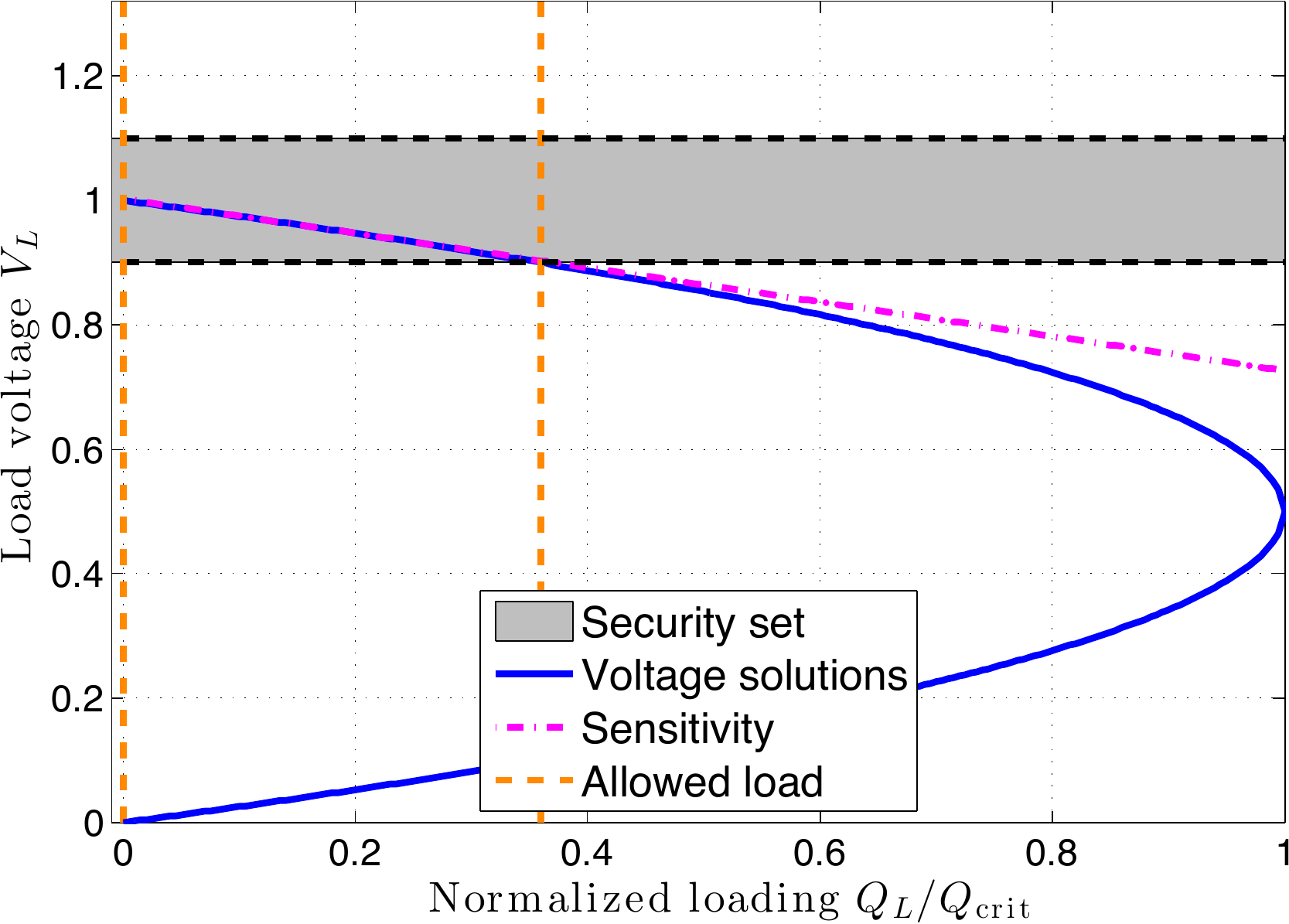}}
\quad
\subfloat[Presence of shunt.\label{subfig:shunts}]{\includegraphics[width=0.4\textwidth]{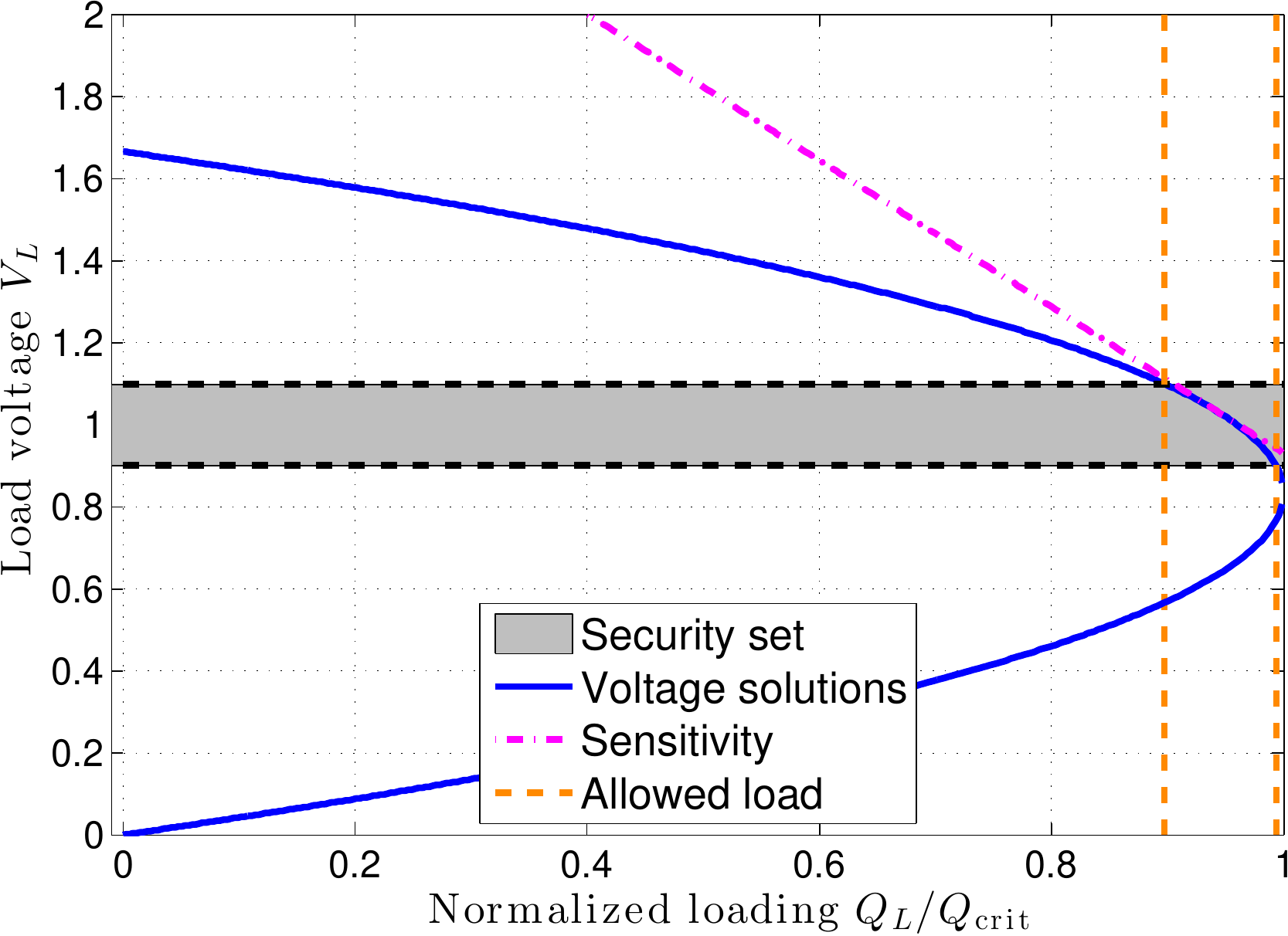}}
\caption{\small Two-buses case: panel (a) plots the network scheme where $V_N=1\ \mathrm{[p.u.]}$, $B=4\ \mathrm{S}$. Panels (b)--(c) plot the QV \emph{nose curve} for two different configurations: (a) Absence of shunt capacitor, $B_{\rm shunt}=0$. (b) Presence of shunt capacitor, $B_{\rm shunt}=2.4\ \mathrm{S}$. 
}
\label{fig:2_node_case_example}
\end{figure*}

\subsection{Decoupled Reactive Power Flow \& Critical Load Matrix}\label{subsec:decoupled_rpf}
Under normal operating conditions, the high-voltage operating point is characterized by small voltages angle differences \cite{JM-JWB-JRB:08} which are treated as parameters \cite{JT-DS-MIS:86} or considered as negligible \cite{RK-FFW:84}. We formalized this statement with the following 

\smallskip

\begin{assumption}[Decoupling Assumption]\label{ass:decoupling}
In steady-state operating conditions, for $\delta\in[0,\pi/2[$, the voltage angle differences are constant and such that
$
|\theta_h - \theta_k | \leq \delta\, ,\ \forall\ (h,k) \in\mathcal{E}.
$
\end{assumption}

\smallskip

Note that, under Assumption \ref{ass:decoupling}, from the form of Eq.\eqref{eq:reactive_pfes}, it is possible to define an \emph{effective susceptance matrix} by embedding the power angle terms into the original line susceptances.
Under the decoupling Assumption \ref{ass:decoupling} the Reactive Power Flow Equations (RPFEs) \eqref{eq:reactive_pfes} simplifies in vector notation to 
\begin{equation}\label{eq:RPFEs}
q(V) = -[V]BV\,.
\end{equation} 
We now take into account the models and the partition introduced in Sections~\ref{subsec:model} and, accordingly we partition the vectors of voltage magnitudes and reactive power injections as
$
V = \begin{bmatrix}
V_L^T & V_G^T 
\end{bmatrix}^T
$,
$
Q = \begin{bmatrix}
Q_L^T & Q_G^T 
\end{bmatrix}^T
$.
%
Combining the power flow \eqref{eq:RPFEs}, the loads model and the partitioning \eqref{eq:Bmatrix}, the power balance $Q_h = q_h(V)$ at each load $h \in \mathcal{V}_L$ can be written as
%
\begin{equation}\label{eq:LoadRPFEs1}
Q_L = -[V_L]\left(B_{LL}V_L + B_{LG}V_G\right)\, .
\end{equation}
We define the \emph{open-circuit voltages} $V_L^*$ as
\begin{equation}\label{eq:VLstar}
V_L^* := -B_{LL}^{-1}B_{LG}V_G\, ,
\end{equation}
which are well defined under Assumption \ref{ass:Mmatrix}. Physically, the open-circuit voltages \eqref{eq:VLstar} are the voltages one would measure at the load buses for zero reactive power demands $Q_L = \vzeros[n_{\ell}]$. With this notation, the RPFEs \eqref{eq:LoadRPFEs1} can be written as
\begin{equation}\label{eq:LoadRPFEs}
Q_L = -[V_L]B_{LL}(V_L - V_L^*)\ .
\end{equation}
Once \eqref{eq:LoadRPFEs} is solved for an operating point $V_L$, the reactive power injections at generators buses $\mathcal{V}_G$ are uniquely determined by substituting the operating point into the final $n_g$ equations in \eqref{eq:RPFEs}.
We define one more useful quantity. 

\smallskip

\begin{definition}[Critical Load Matrix]\label{def:critical_load}
Given the matrix $B_{LL}$ and the open-circuit profile $V_L^*$ as defined in \eqref{eq:VLstar}, the critical load matrix $\Qcrit$ is defined as
\begin{equation}\label{eq:CriticalMatrix}
\Qcrit := \frac{1}{4}[V_L^*]B_{LL}[V_L^*].
\end{equation}
\end{definition}

The critical load matrix $Q_{\rm crit}$ concisely combines the network structure, generator voltages, shunts, and the relative locations of generation and load. In particular, it will help us to formulate the optimal voltage support problem to follow.
%
Finally, it is convenient to rewrite \eqref{eq:LoadRPFEs} in a normalized set of variables. Using the open-circuit voltages $V_L^*$ defined in \eqref{eq:VLstar}, we denote the vector of \emph{normalized voltages} as 
\begin{equation}\label{eq:normalized_coordinate}
v := [V_L^*]^{-1}V_L\,.
\end{equation}
Note that if the open-circuit profile is flat ($V_L^* = \alpha \vones[]$ for some $\alpha > 0$), then $v$ is simply the standard vector of per unit voltages. In general, however, due to inhomogeneous  generators voltage set points and the presence of shunt compensation, $V_L^*$ is not flat and the scalings in \eqref{eq:normalized_coordinate} are non-uniform.
Substituting $V_L = [V_L^*]v$ into the RPFEs \eqref{eq:LoadRPFEs} and using \eqref{eq:CriticalMatrix}, \eqref{eq:LoadRPFEs} takes the simple form
\begin{equation}\label{eq:LoadRPFEnormalize}
Q_L = -4[v]\Qcrit(v-\1).
\end{equation}

\section{Formulation of Optimal Voltage Support Problem}\label{sec:problem}
In this section we present our novel problem formulation.
First, we present a common operational requirement highlighting its possible inadequacy to capture the safe operation of the grid.
Then, we present a novel metric to measure the stress induced on the network by the load demand, representing our objective function. Finally, we present a linearization which leads to a convex reformulation of the optimization problem.

\subsection{Security Constraints}\label{subsec:security_constraint}
A common operational requirement is that the load buses voltage magnitudes must lie within a predefined percentage 
deviation, typically $5\%$, from a reference voltage. 
This tight clustering of voltages is due to the following reasons:
\begin{enumerate}[(i)]
\item loads and some system components are designed to operate with a voltage in a narrow region around the network base voltage; 
\item a flat voltage profile minimizes current flows 
and, consequently, minimizes resistive power losses;
\item a flat profile usually reduces the sensitivity of the voltage profile with respect to load changes (see Example \ref{ex:secure_set});
\item most importantly, by conventional wisdom (see typical nose curve studies \cite{TVC-CV:98,ARB-VV:06}) a flat voltage profile indicates that the network is safe from voltage collapse. 
\end{enumerate}
%
We formalize this requirement by defining the \emph{secure set}.

\smallskip

\begin{definition}[Secure set]\label{def:secure_set}
Given a reference voltage $V_N~\in~\mathbb{R}_{>0}$, a percentage deviation $\alpha>0$ and $V_L^*$ as in \eqref{eq:VLstar}, the \textit{secure set} $\secure_set$ is defined as
\begin{equation}\label{eq:secure_set}
\secure_set := \left\{ v\in\mathbb{R}^{n_\ell}\ \Big|\ \frac{ \left\| [V_L^*]v - V_N\1 \right\|_{\infty}}{V_N}\leq\alpha \right\}\ .
\end{equation}
\end{definition}

\smallskip

Hence, if $v \in \secure_set$ is a solution to \eqref{eq:LoadRPFEnormalize}, then all voltages lie within $\alpha$ percent of the nominal voltage $V_N$. 
While this represents a baseline operational requirement, under some circumstances it may not be sufficient to ensure safe grid operation. We present a simple example highlighting this fact.

\smallskip

\begin{example}[Security requirement inadequacy]\label{ex:secure_set}
Consider the simple two-buses case study consisting of a load connected to a source at voltage $V_N=1$, as illustrated in Figure~\ref{subfig:2node_case}.
For the case where $B_{\rm shunt} = 0$, Figure \ref{subfig:no_shunts} plots the \emph{nose curve}, i.e., the locus of solutions to \eqref{eq:LoadRPFEs} (blue solid blue) as $Q_L$ is varied from $0$ to $Q_{\rm crit}$. Note that for a chosen $Q_L$, there may be two, one, or zero feasible solutions of \eqref{eq:LoadRPFEs}. The secure set is shown as a shaded area between two dashed black lines. Also shown are the loading limits which ensure the high-voltage solution to lie in the secure set (dashed orange), and the tangent line to the nose curve at the mid-point between the dashed orange lines (dashed magenta). This tangent line captures the sensitivity of the load voltage to changes in reactive power demand. From Figure \ref{subfig:no_shunts}, note that if $Q_L$ is too large,  the operating point does not lie within $\secure_set$. A standard policy is then to support the voltage level by adjusting the shunt compensation, i.e., by increasing $B_{\rm shunt}$. When $B_{\rm shunt} = 0$, Figure~\ref{subfig:no_shunts} demonstrates that the security requirement $v\in\secure_set$ guarantees a \tq{safe} distance to collapse, represented by the nose of the blue curve. Moreover, the sensitivity of the voltage to changes in load is small, meaning that relatively large changes in loading do not translate into large voltage changes.
Conversely, in Figure~\ref{subfig:shunts} $B_{\rm shunt} \neq 0$, and the security requirement $v \in \secure_set$ is \tq{dangerously} close to the nose of the curve. 
%
%
Finally, the sensitivity line is steeper meaning that small changes in the load cause relatively big changes in the voltage. This affects the robustness of the network to small load changes.
%
\oprocend
\end{example}

\smallskip 

%
The previous analysis highlights that the security requirement $v\in\secure_set$ alone could be insufficient. Note that, to operate the grid in the point farthest from voltage collapse and to ultimately maximize the stability and robustness margins, a simple intuition is that of minimizing the distance of the operating point from the open-circuit solution $V_L^*$, represented by the left-most point on the blue curve \textemdash{} constrained to the fact that the operating point must belong to $\secure_set$. As final remark, note that in general $V_L^*$ does not belong to $\secure_set$. Thus in general, distance-to-collapse minimization and voltage compensation do not coincide.

\subsection{Network Stress Measure and Stress Minimization Problem}\label{subsec:stress_measure}
Based on the insights given by Example \ref{ex:secure_set}, we define the following measure quantifying the distance to 
collapse.

 
\smallskip

\begin{definition}[Network Stress Measure]\label{def:stress_measure}
Consider the RPFEs \eqref{eq:LoadRPFEnormalize} in the normalized voltages $v\in\R^{n_\ell}$. The \emph{network stress measure} induced by the load is defined as
\begin{equation}\label{eq:stress_measure_full}
J_{\rm stress}(v) := \left\| v - \1 \right\|_{\infty}\, .
\end{equation}
\end{definition}

\smallskip

Definition \ref{def:stress_measure} is based on the intuition that the open-circuit profile $V_L^*$ is the network's natural operating point in absence of loading, i.e., under \tq{no stress}. Conversely, when the network works close to the nose tip, i.e., the farthest point from $V_L^*$ then, this is a \tq{high-stress} scenario. In this sense, the stress function \eqref{eq:stress_measure_full} quantifies the loading on the network conveniently expressed in the normalized profile $v = [V_L^*]^{-1}V_L$.

In the following, we assume that a certain number of load buses can be equipped with additional controlled devices, e.g., synchronous condensers \cite{GA:08}. 
We assume these devices can provide a controllable amount of reactive power support, and in the following we model them as controllable sources of reactive power $q_h$, subject to upper and lower operational bounds. Specifically, the RPFEs  \eqref{eq:LoadRPFEnormalize} are modified as 
\begin{equation}\label{eq:LoadRPFEnormalize_controlled}
Q_L + q = -4[v]Q_{\rm crit}(v-\1)\ ,
\end{equation}
where $q\in\R^{n_\ell}$ is such that $q_{\min}\leq q\leq q_{\max}$ and $q_{\min},q_{\max}\in\R^{n_\ell}$ are vectors representing the injection capacity constraints. If load bus $h \in \mathcal{V}_L$ is not equipped with a compensator, we set $q_{\min , h} = q_{\max , h} = 0$.

We now formulate our optimization problem of interest, which we refer to as the \textit{Stress Minimization} problem.

\smallskip

\begin{problem}[Stress Minimization]\label{prob:stress_minimization}
Given $Q_L$, $\Qcrit$ as in \eqref{eq:CriticalMatrix} and the capacity limits $q_{\min}$, $q_{\max}$, find $q$ and $v$ such that 
%
\begin{eqnarray}\label{eq:stress_minimization}
& \underset{q\in\mathbb{R}^{n_\ell}}{\mathrm{minimize}}	& 	J_{\rm stress}(v)  ,\\
& \mathrm{subject\ to}	& 	\begin{cases}
						v\in\secure_set\, ,\\
						q_{\rm min} \leq q \leq q_{\rm max}\, ,\\
					Q_L + q = -4[v]Q_{\rm crit}(v-\1)
					\, . 
			\end{cases}	\nonumber
\end{eqnarray}
\end{problem} 

\smallskip

The main idea behind Problem \ref{prob:stress_minimization} is that minimizing $J_{\rm stress}$ keeps the operating point away from the tip of the nose curve, i.e., the collapse point. 
The standard security requirement $v~\in~\secure_set$ is imposed as a hard constraint.

%
Since $v$ is related to $q$ through the quadratic equality constraints \eqref{eq:LoadRPFEnormalize_controlled}, Problem \ref{prob:stress_minimization} is nonlinear and non-convex.
In the following, we convexify this problem through the use of a power flow linearization.

\begin{figure*}[t]
\centering
\subfloat[$\gamma = 0$]{\includegraphics[width = 0.33\textwidth]{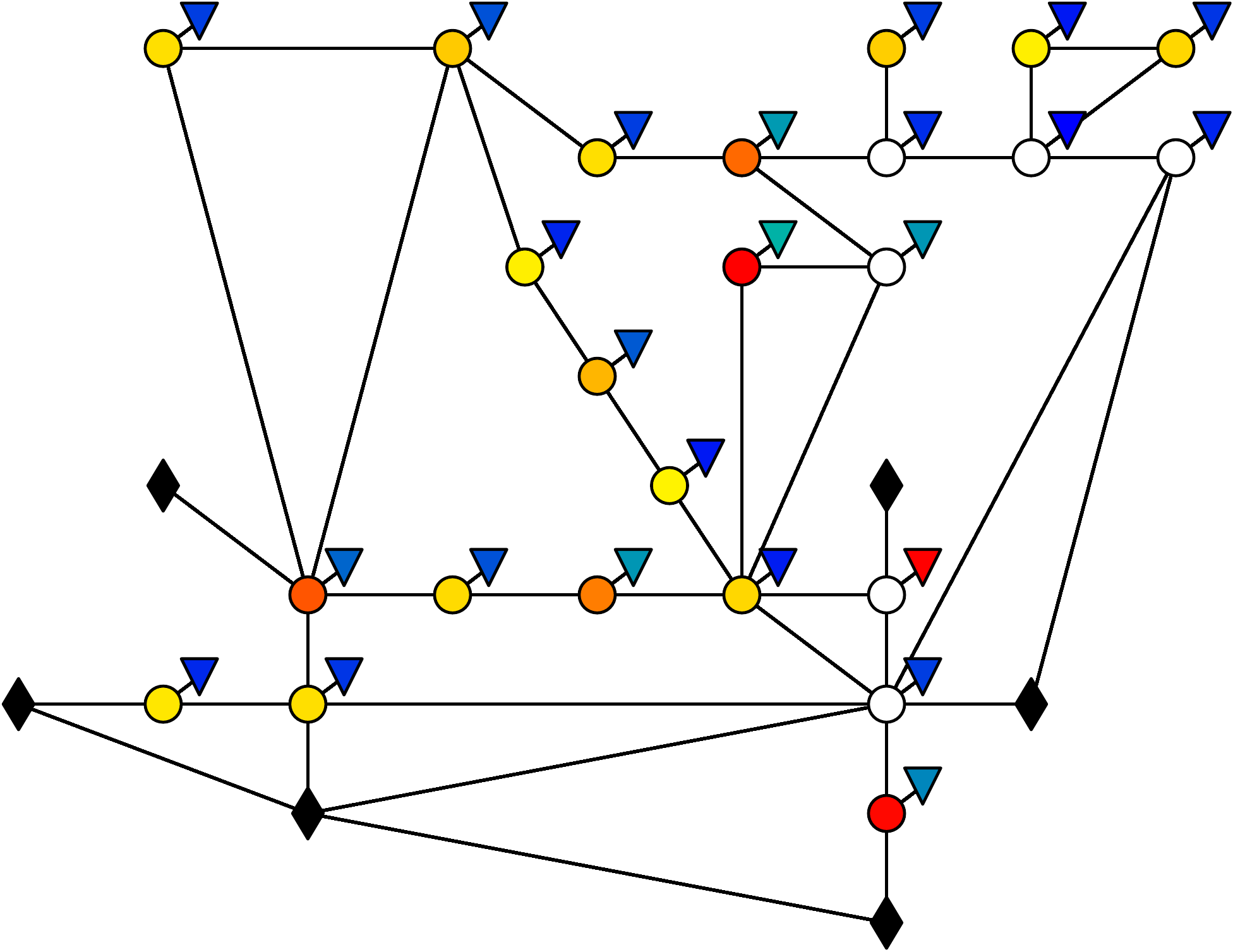}\label{subfig:gamma0}}
\subfloat[$\gamma = 4\times10^{-4}$]{\includegraphics[width = 0.33\textwidth]{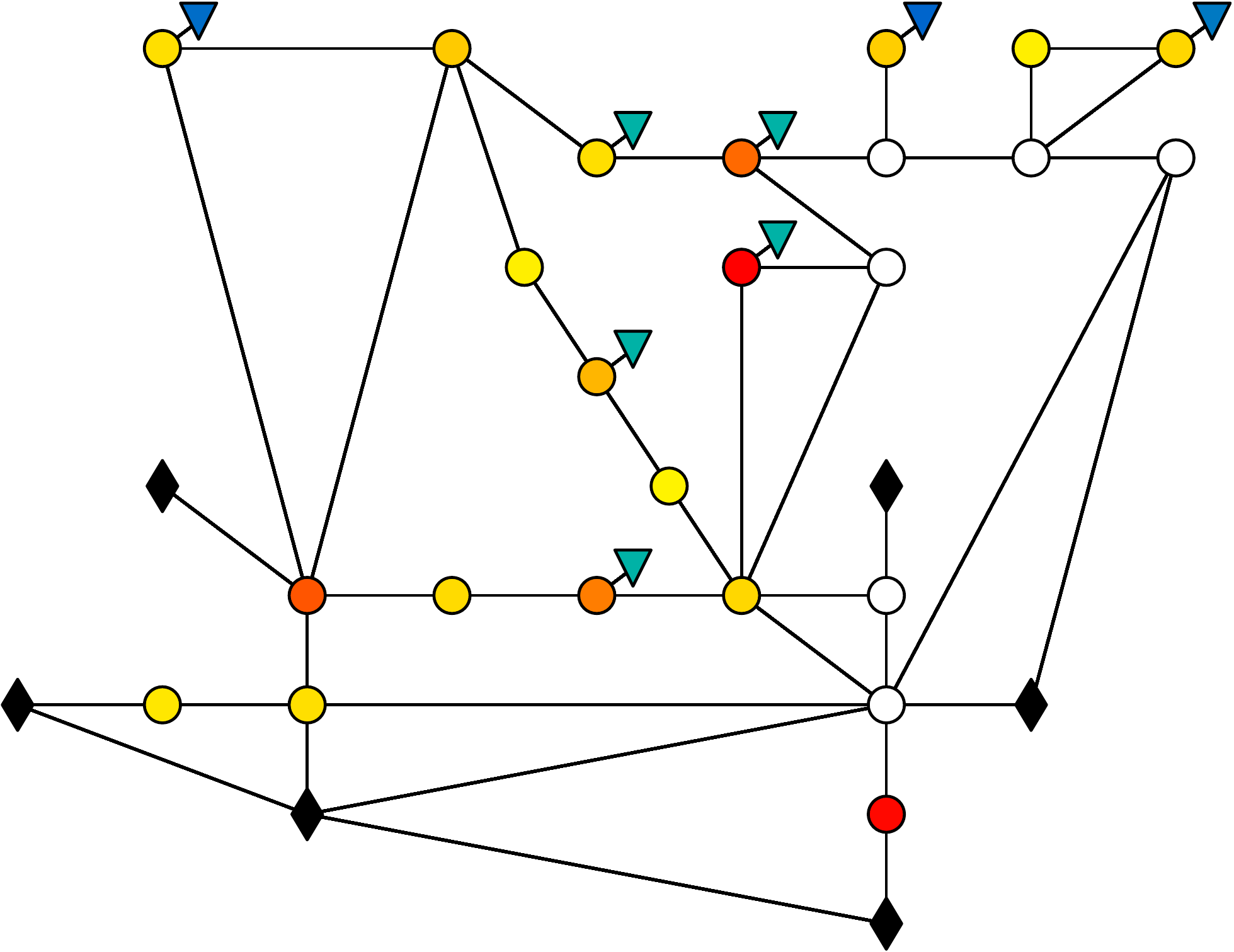}\label{subfig:gamma4}}
\subfloat[$\gamma = 8\times10^{-4}$]{\includegraphics[width = 0.33\textwidth]{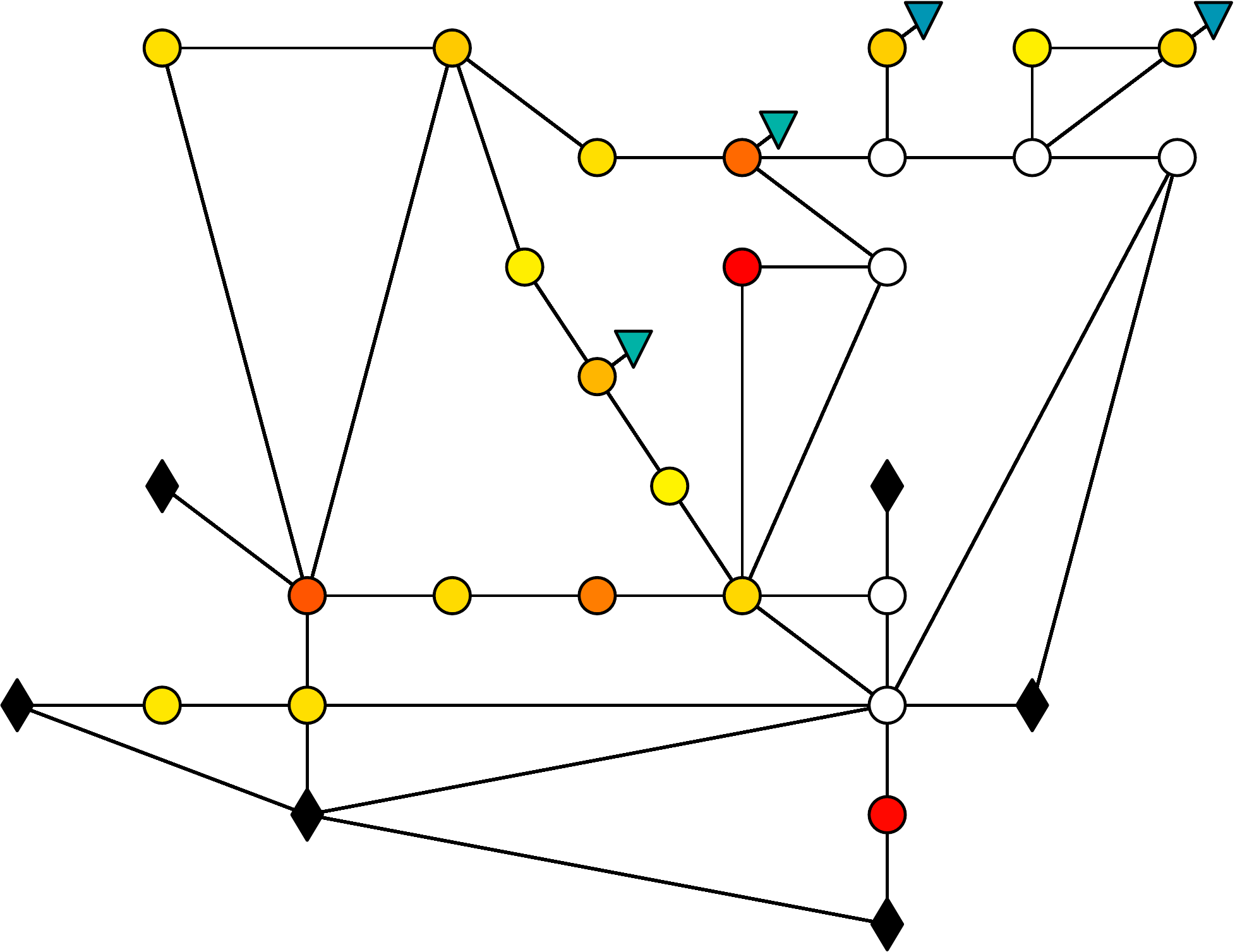}\label{subfig:gamma8}}
\caption{\small Placement scheme  of the compensation devices for different values of $\gamma$. Black diamonds $\vardiamond$ represent generators, circles $\bigcirc$ represent loads, and triangles $\bigtriangledown$ represent reactive compensators. Color scheme: red scale $\protect\mycbox{golden} - \protect\mycbox{red}$ represents increasing power absorption; blue scale $\protect\mycbox{blue} - \protect\mycbox{forestgreen}$ represents increasing power injection.}
\label{fig:GridGraphicVisualization}
\end{figure*}

\subsection{Linear Approximation and Convexification}\label{subsec:linearization}

We now introduce a suitable linearization for Problem \ref{prob:stress_minimization} which had been first presented~in~\cite{BG-JWSP-FD-SZ-FB:13zb}.
From \eqref{eq:LoadRPFEnormalize_controlled}, assuming $\|Q_L+q\|\sim 0$, we expect the normalized profile \eqref{eq:normalized_coordinate} to be $v\simeq \1$ which would be the exact high voltage solution corresponding to $Q_L+q=0$.
Linearizing the RPFEs \eqref{eq:LoadRPFEnormalize_controlled} around $v = \1$, to first order, the solution is given by
\begin{equation}\label{eq:solution_first_order_approx_controlled}
\widehat{v} = \1 -\frac{1}{4}\Qcrit^{-1} (Q_L+q)\ .
\end{equation}
That is, to first order the solution of \eqref{eq:LoadRPFEnormalize_controlled} is given by a uniform component plus a deviation which is linear in the reactive injections.
%
Using \eqref{eq:solution_first_order_approx_controlled}, the cost \eqref{eq:stress_measure_full} is approximated by
\begin{equation}\label{eq:stress_measure_approx}
J_{\rm stress}(v) = \left\| v - \1 \right\|_{\infty} \propto \left\| \Qcrit^{-1}(Q_L+q) \right\|_{\infty}.
\end{equation}
Note that the approximated cost function \eqref{eq:stress_measure_approx} is convex in the reactive power injections $q$. 
%
%
%
%
%
By exploiting \eqref{eq:solution_first_order_approx_controlled} and thanks to some algebraic manipulations, it is possible to see that the security requirement $\widehat{v}\in\secure_set$ holds if and only if
\begin{equation}\label{eq:sec_constr}
V_N(1 - \alpha)[V_L^*]^{-1}\1 \leq \widehat{v} \leq V_N(1 + \alpha)[V_L^*]^{-1}\1.
\end{equation}
Substituting for $\widehat{v}$ from \eqref{eq:solution_first_order_approx_controlled}, \eqref{eq:sec_constr} is equivalent to 
\begin{equation}\label{eq:security_constraints}
\xi_{\min} \leq -Q_{\rm crit}^{-1}q \leq \xi_{\max}\ ,
\end{equation}
where 
\begin{subequations}
\begin{equation}\label{eq:security_threshold_low}
\xi_{\min} := 4\left(V_N \left( 1-\alpha \right) [V_L^*]^{-1}\1 - \1 \right) + \Qcrit^{-1}Q_L\  ,
\end{equation}
\begin{equation}\label{eq:security_threshold_up}
\xi_{\max}\ := 4\left(V_N \left( 1+\alpha \right) [V_L^*]^{-1}\1 - \1 \right) + \Qcrit^{-1}Q_L\ .
\end{equation}
\end{subequations}
%
Thus, the security constraints $\widehat{v} \in \secure_set$ have been converted into linear inequality constraints  on the decision variables $q$. We now present the convexified version of Problem \ref{prob:stress_minimization}.

\smallskip

\begin{problem}[Convex Stress Minimization]\label{prob:stress_minimization_approx}
Consider the RPFEs \eqref{eq:LoadRPFEnormalize_controlled}. Let $\Qcrit$ be as in \eqref{eq:CriticalMatrix} and $q_{\rm min}$, $q_{\rm max}$  be vectors representing the injection capacity limits. Finally, define $\xi_{\min}$ and $\xi_{\max}$ as in \eqref{eq:security_threshold_low}--\eqref{eq:security_threshold_up}, respectively.
%
Then the goal is to
\begin{eqnarray}\label{eq:stress_minimization_approx}
&	\underset{q\in\mathbb{R}^{n_\ell}}{\mathrm{minimize}}	& 	\left\| \Qcrit^{-1}\left(Q_L + q\right)\right\|_{\infty},\\
&	\mathrm{subject\ to}		& 	\begin{cases}\label{eq:convex_constraints}
							\xi_{\min} \leq -\Qcrit^{-1}q \leq \xi_{\max}\,,\\
							q_{\rm min} \leq q \leq q_{\rm max}\,.
							\end{cases}
\end{eqnarray}
\end{problem} 

\smallskip

\begin{remark}[On the stress measure]\label{rem:stress_measure}
Aside from the linearization-based derivation in this subsection, the measure \eqref{eq:stress_minimization_approx} is inspired by recent results \cite{JWSP-FD-FB:14c} on the solvability of the decoupled reactive power flow equations \eqref{eq:LoadRPFEs},
where it has been shown that $\|\Qcrit^{-1}Q_L\|_{\infty}$ represents a proper distance-to-collapse measure. Indeed, if
$
\|\Qcrit^{-1}Q_L\|_{\infty} < 1\,,
$
the non linear \eqref{eq:LoadRPFEs} has a unique high-voltage solution safe from collapse.
\oprocend
\end{remark}

\smallskip

Observe that in Problem \ref{prob:stress_minimization_approx} the cost
\eqref{eq:stress_minimization_approx} and the constraints
\eqref{eq:convex_constraints} are convex in the decision variables.
Moreover, Problem~\ref{prob:stress_minimization_approx} can be written as a linear program 
and can therefore be efficiently solved via convex optimization.

Before presenting some performance and simulations of the stress minimization procedure, notice that both Problems \ref{prob:stress_minimization} and \ref{prob:stress_minimization_approx} are offline centralized procedures which, as suggested by the formulation in Section \ref{subsec:stress_measure}, assume that either the full set of load buses or only an a priori assigned subset of them are equipped with controllable devices. The first scenario is impractical and economically unfeasible in large networks due to the large number of devices needed. The second scenario could likely lead to a sub-optimal allocation of resources if no specific allocation policies are used. In the following subsection, we refine Problem~\ref{prob:stress_minimization_approx} to simultaneously solve for the planning problem of allocating the resources along with the system-level stress minimization problem.  

\subsection{The Planning Problem: Sparse Stress Minimization}\label{subsec:sparsity_promotion}
%
Here we propose a modification of Problem \ref{prob:stress_minimization_approx} to find a desired trade-off between the number of actuators and the minimization of the stress cost. In order to accomplish this task, we propose a \emph{sparsity-promoting} approach where, by tuning an additional parameter, the user is able to control the sparsity of the solution. In this way, we simultaneously solve the system-level stress minimization problem as well as the planning problem of allocating a finite number of resources.

%
The cardinality function, $\mathrm{card}(\cdot)$, is a natural choice to account for the number of devices. However, it is discontinuous and non-convex. 
A convex approximation of $\mathrm{card}(q)$ is the \textit{re-weighted} $\ell_1$-norm \cite{SB-LV:04} 
\begin{equation}\label{eq:re-weighted_norm}
\left\| [w(q)]q \right\|_1 = \sum_{h = 1}^{n_\ell} w_h(q_h)q_h\, ,\ \ \ \ w_h(q_h) := \frac{1}{|q_h| + \epsilon}\, ,
\end{equation}
where $0 < \epsilon \ll 1$. Adding equation \eqref{eq:re-weighted_norm} to the cost function \eqref{eq:stress_minimization_approx}, it is possible to formulate the following problem which we refer to as the \textit{Sparse Stress Minimization} problem.

\smallskip

\begin{problem}[Sparse Stress Minimization]\label{prob:sparse_stress_minimization_approx}
Consider the same set-up as in the Convex Stress Minimization of  Problem~\ref{prob:stress_minimization_approx}.
%
Then the goal is to
\begin{eqnarray}\label{eq:sparse_stress_minimization_approx}
&	\underset{q\in\mathbb{R}^{n_\ell}}{\mathrm{minimize}}	& 	\left\| \Qcrit^{-1}\left(Q_L + q\right)\right\|_{\infty} + \gamma\left\| [w(q)]q \right\|_1,\\
&	\mathrm{subject\ to}		& 	\begin{cases}
							\xi_{\min} \leq -\Qcrit^{-1}q \leq \xi_{\max},\\
							q_{\rm min} \leq q \leq q_{\rm max}.
							\end{cases}							
							\nonumber							
\end{eqnarray}
\end{problem}

\smallskip

The parameter $\gamma$ in the cost function \eqref{eq:sparse_stress_minimization_approx} can be used to promote sparsity of the solution $q$, and thereby minimize the number of required actuators. Obviously for $\gamma = 0$, Problem~\ref{prob:sparse_stress_minimization_approx} reduces to Problem~\ref{prob:stress_minimization_approx}. By increasing the value of $\gamma$ the user can force the solver to lean towards a more sparse solution. This automatically compels the solver to optimally allocate the resources in order to find the best trade-off between sparsity and system-level stress minimization.


\subsection{Simulation: Planning Problem and Offline Optimization}\label{subsec:simulation_sparsity} 
We now present a case study to show the effectiveness of planning and
the offline optimization procedure proposed. The simulations refer to
Problem \ref{prob:sparse_stress_minimization_approx} and are
implemented in MATLAB and CVX ~\cite{MG-SB:11-cvx}. The plotted
voltage profiles refer to the linearized solution
\eqref{eq:solution_first_order_approx_controlled} of the decoupled
RPFEs \eqref{eq:LoadRPFEnormalize_controlled}. The test-bed consists
of:
\begin{itemize}
\item IEEE 30 bus transmission grid \cite{RDZ-CEMS-RJT:11};
\item a reference voltage $V_N = 1$ [p.u.];
\item a voltage deviation limit $\alpha=5\%$;
\item capacity limits $\{q_{\rm min},q_{\max}\} = \{-0.5,0.5\}\times\|Q_L\|_{\infty}\1$.
\end{itemize}
%
\begin{figure}[t]
\centering
\includegraphics[width = \columnwidth]{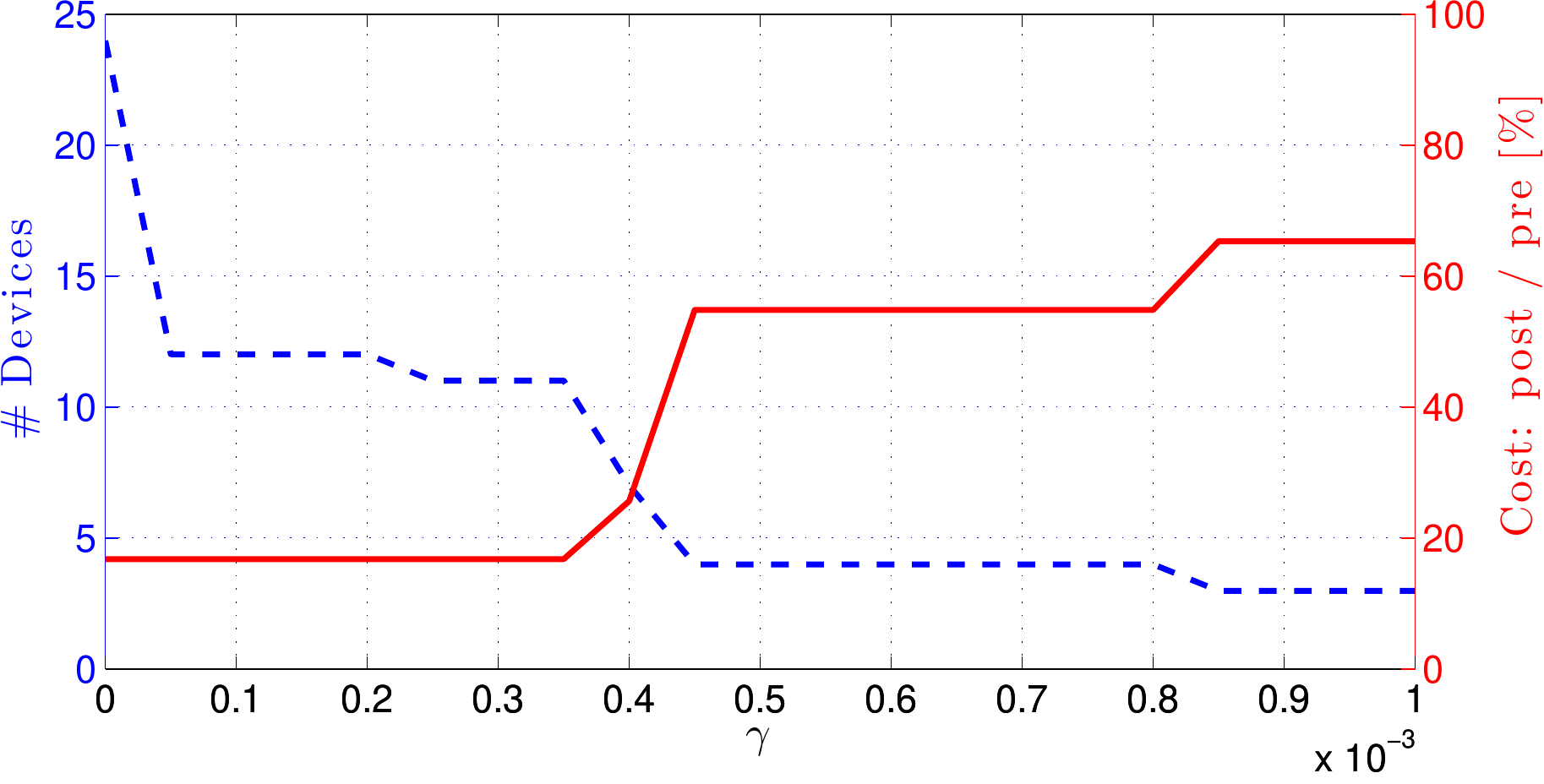}
\caption{\small Behavior of the number of devices (left-blue axis) and of the ratio of the cost value after and before the optimization (right-red axis) as function of the sparsity parameter $\gamma$.}
\label{fig:Behavior_NumberDevice_CostImprovement}
\end{figure} 
%

\begin{figure}
\centering
\includegraphics[width = \columnwidth]{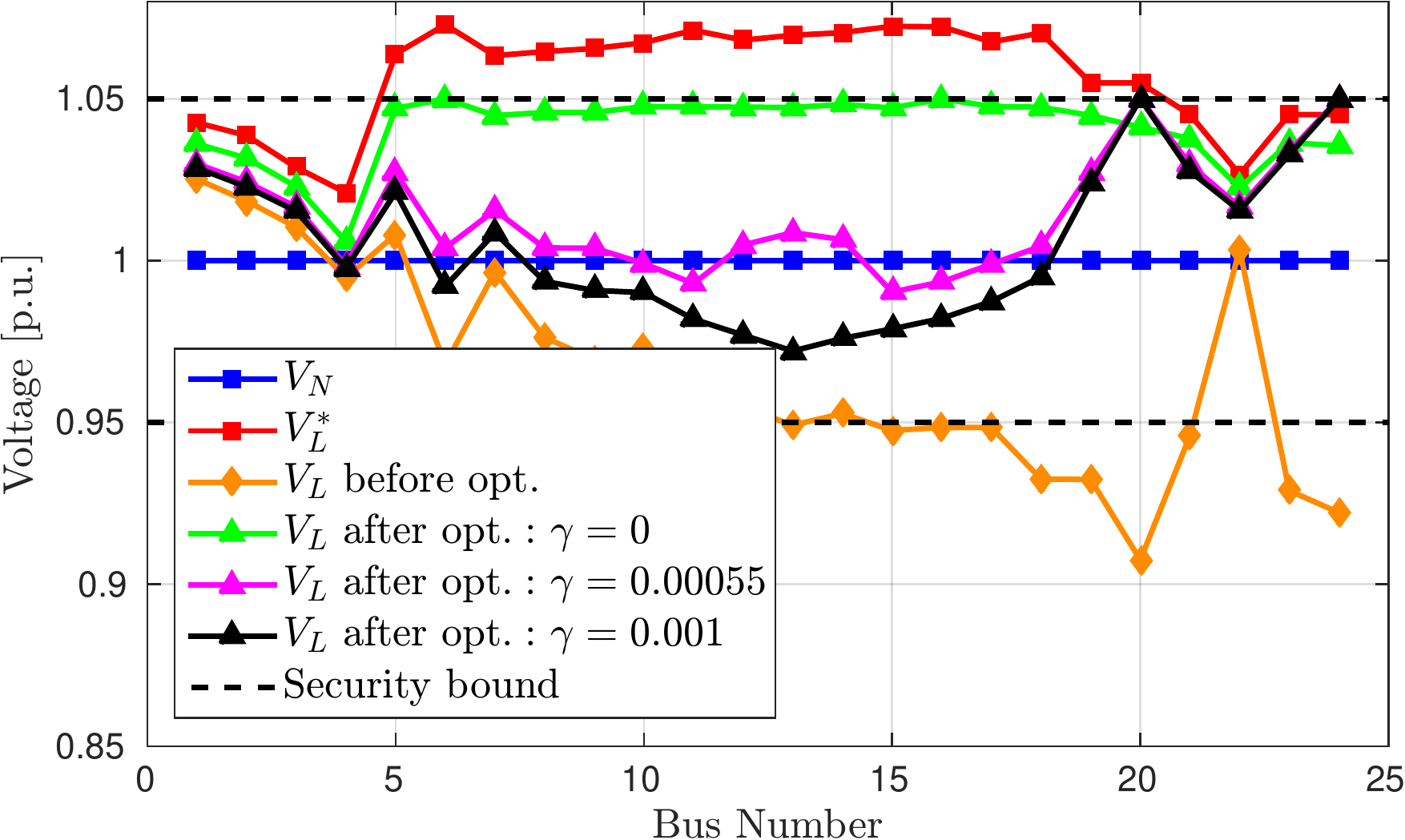}
\caption{\small Voltage profiles for different values of $\gamma$.}
\label{fig:profile_before_after}
\end{figure}

Figures \ref{subfig:gamma0}--\ref{subfig:gamma4}--\ref{subfig:gamma8} 
illustrates the placement of actuators
for increasing values of $\gamma$ ($\gamma = 0$, $\gamma = 4\times 10^{-4}$ and $\gamma = 8\times 10^{-4}$, respectively). Compensators placed by the optimization problem are indicated with a triangle.
The color scheme for loads $Q_L$ and compensators $q$ is as follows:
\begin{itemize}
\item reactive injections, i.e., positive values, are plotted in blue-scale ($\protect\mycbox{blue} - \protect\mycbox{forestgreen}$): the lighter the blue, the bigger the injection absolute value;
\item reactive consumptions/absorptions, i.e., negative values, are plotted in red-scale ($\protect\mycbox{golden} - \protect\mycbox{red}$): the darker the red, the bigger the absolute value of the consumption;
\item white node means zero injections/consumptions.
\end{itemize}
First of all, for $\gamma=0$ the solver places compensation everywhere. Moreover it can be seen that one compensator is red colored, meaning that it effectively absorbs reactive power. 
This occurs due to large reactive power injections at neighboring buses, which drive up voltage values across the network -- additional reactive power must be absorbed to lower specific voltages and meet the security constraints.
%
%
Sparsity promotion takes place for increasing $\gamma$ with the solver placing compensators only 
where the heaviest loading occurs.
Figure~\ref{fig:Behavior_NumberDevice_CostImprovement} shows, as a function of $\gamma$, the number of devices placed (left-blue axis) and the ratio between the value of the cost \eqref{eq:stress_minimization_approx} after a polishing step, i.e., obtained as solution of Problem \ref{prob:stress_minimization_approx} given the placement obtained solving Problem \ref{prob:sparse_stress_minimization_approx}, over the value of the cost before the optimization (right-red axis), i.e., 
$$
\frac{\|\Qcrit^{-1}(Q_L + q)\|_\infty}{\|\Qcrit^{-1}Q_L\|_\infty}\, .
$$
It can be seen how, for increasing $\gamma$ the final value of the cost increases since a smaller number of controllable units are less able to compensate the voltage profile. However, as can be seen from the first part of the plot, by using only 11 controllers we achieve the same performance as of using 24 compensators. This not only highlights the redundancy of using 24 compensators but that the optimal placement is necessary to achieve the same level of performance. 
Figures~\ref{fig:profile_before_after} shows the linearized profiles of $V_L$ before and after the optimization for different $\gamma$. It can be seen that for increasing $\gamma$ the profile $V_L$ is less compensated, i.e., it is farther from the $V_L^*$ profile.
Finally, as already stressed, note that stress minimization and classical voltage compensation do not coincide. Indeed, $V_L^*$, in general, does not belong to $\secure_set$, identified by the black dashed lines.
%

\section{Distributed Online Stress Minimization through Feedback Control}\label{sec:distributed_algorithm}
The previous methods of Sections \ref{subsec:linearization} and \ref{subsec:sparsity_promotion} are suitable only for offline optimization and planning. In this section we assume the planning problem has been solved offline, and develop a dual-ascent algorithm for Problem \ref{prob:stress_minimization_approx} which may be implemented online as a distributed feedback controller. This is motivated by different reasons among which it is worth mentioning that:
\begin{itemize}
\item utilities could prefer not to share information with a central operator because of privacy reasons;
\item an online implementation can be naturally exploited as a distributed feedback controller in presence of time-varying loads, to reject disturbances, to increase the system robustness, and to track the optimal solution.
\end{itemize}
In the following, we assume that \tq{smart agents} are embedded at all the grid's load buses. These are characterized by mild communication and computational capabilities. Moreover, they can communicate according to a communication graph which is designed to coincide with the electrical network. It is worth mentioning that, as will be clear later, even the load buses not equipped with a controllable compensator are required to share \tq{smartness} capabilities.

In the current formulation of Problem \ref{prob:stress_minimization_approx}, the presence of the dense matrix $\Qcritinv$ in both the cost $J_{\rm stress}$ and the constraints \eqref{eq:security_constraints} compromises the possibility to solve the stress minimization problem in a distributed fashion. However, the matrix $\Qcrit$ is sparse and the graph induced by its sparsity pattern coincides with the topology of the grid which connects the load buses. We take advantage of this structure to develop a distributed algorithm to solve Problem \ref{prob:stress_minimization_approx}.

Whereas the formulation of Problem \eqref{prob:stress_minimization_approx} expressed the stress minimization compactly in injection coordinates $q$, now we derive the equivalent formulation in voltage coordinates to leverage on the sparsity of $\Qcrit$.
%
%
%
We start our analysis defining the \emph{deviation variable} $x$ as
\begin{equation}\label{eq:voltage_deviation}
x := -\Qcritinv(Q_L + q)\, ,
\end{equation}
which represents the linear deviations of the voltages $\widehat{v}$ as defined in \eqref{eq:solution_first_order_approx_controlled} due to the overall reactive injection.
Similar to what done in Section \ref{subsec:linearization}, from the definition of set $\secure_set$ it is possible to obtain the security constraints expressed in the $x$ coordinates. These are equal to
$$ x_{\min} \leq x \leq x_{\max}\, ,$$
where 
\begin{subequations}
$$
x_{\min} := 4\left(V_N \left( 1-\alpha \right) [V_L^*]^{-1}\1 - \1 \right)\,,
$$
$$
x_{\max} := 4\left(V_N \left( 1+\alpha \right) [V_L^*]^{-1}\1 - \1 \right)\,.
$$
\end{subequations}
From the definition of $x$ it is clear that
\begin{equation}\label{eq:control_input_distributed}
q = -(\Qcrit x - Q_L)\, .
\end{equation}
Since the matrix $\Qcrit$ is characterized by a sparsity pattern equivalent to that induced by the electric graph connecting the loads, the desired control inputs can be computed by means of a local exchange of information, namely the $x_i$ variables among electric neighbors. Additionally, from \eqref{eq:control_input_distributed} it is easy to impose the capacity constraints, i.e., $q_{\min}\leq q\leq q_{\max}$. The Problem \ref{prob:stress_minimization_approx} is then equivalent to

\smallskip

\begin{problem}[Online Stress Minimization]\label{prob:distributed_stress_minimization_approx}
\begin{eqnarray}\label{eq:distributed_stress_minimization_approx}
&	\underset{x\in\mathbb{R}^{n_\ell}}{\mathrm{minimize}}	& 	\left\| x \right\|_{\infty},\\
&	\mathrm{subject\ to} 	& 	\begin{cases}
							x_{\min} \leq x \leq x_{\max},\\
							q_{\rm min} \leq -(\Qcrit x - Q_L) \leq q_{\rm max}.
							\end{cases}							
							\nonumber							
\end{eqnarray}
\end{problem}

\smallskip

Now, we point out three more issues related to Problem \ref{prob:distributed_stress_minimization_approx}: (i) the $\infty$-norm is not everywhere differentiable and thus not suitable for a gradient-based iterative procedure; (ii) computing the cost in \eqref{eq:distributed_stress_minimization_approx} requires knowledge of all $x_i$ variables; (iii) in order to compute the derivative of the maximum function embedded in $\|\cdot\|_{\infty}$, the index where the maximum is attained must be known.
Next, we propose one possible solution to these issues. 

\subsection{A Smooth Decomposable Approximation of $\infty$-norm}\label{subsec:smooth_infty_norm} 
We now present a continuously differentiable approximation for the $\infty$-norm which combines a smooth approximation for the maximum function, the softmax \cite{CMB:06}, and a smooth approximation for the absolute value. This reads as
\begin{eqnarray}\label{eq:smooth_infty}
\tilde{f}_{\alpha,\epsilon}(x)&:=&\mathrm{softmax}_\alpha(|x|^{1+\epsilon})\, ,\ \ 1\ll\alpha\ ,\ 0<\epsilon\ll 1\, , \notag\\
	&=& \frac{1}{\alpha} \log \left( \frac{1}{n_\ell}\sum\nolimits_{i=1}^{n_\ell} \exp\left(\alpha |x_i|^{1+\epsilon}\right) \right)\, .
\end{eqnarray}
The idea behind \eqref{eq:smooth_infty} is to exploit the super-linearity property of the exponential to let the maximum component of the vector $x$ dominate the other components. The exponentiation is exploited to recover differentiability of the absolute value. The approximation $\tilde{f}_{\alpha,\epsilon}(x)$ approximates $\|x\|_{\infty}$ in the following sense; the proof can be found in Appendix \ref{app:proof_lemma_recover_infty_norm}.

\smallskip 

\begin{lemma}[Limit behavior of $\tilde{f}_{\alpha,\epsilon}$]\label{lemma:recover_infty_norm}
Consider $\tilde{f}_{\alpha,\epsilon}$ as in \eqref{eq:smooth_infty}. Then, it holds that
$$
\underset{\substack{\alpha\rightarrow+\infty\\ \epsilon\rightarrow 0^+}}{\lim}\ \tilde{f}_{\alpha,\epsilon}(x) = \|x\|_{\infty}\, .
$$
\end{lemma}

\smallskip 

Now, let us define
\begin{equation}\label{eq:cost_function}
f_{\alpha,\epsilon}(x) := n_\ell\exp(\alpha \tilde{f}_{\alpha,\epsilon}) =  \sum_{i=1}^{n_\ell} \exp(\alpha |x_i|^{1+\epsilon})\,.
\end{equation}
Notice that, due to the monotonicity of the $\log$ function in the softmax, it holds that
\begin{equation}\label{eq:argmin_equiv}
\underset{x\in\R^{n_\ell}}{\mathrm{argmin}}\ \tilde{f}_{\alpha,\epsilon}(x) \equiv 
\underset{x\in\R^{n_\ell}}{\mathrm{argmin}}\ f_{\alpha,\epsilon}(x)\, .
\end{equation}
Thanks to this simplification and in view of the distributed implementation, notice that the $i$-th component of the gradient vector of \eqref{eq:cost_function} is equal to
\begin{equation}\label{eq:cost_derivative}
\frac{\partial f_{\alpha,\epsilon}(x)}{\partial x_i} = \alpha(1+\epsilon)\exp(\alpha |x_i|^{1+\epsilon})|x_i|^\epsilon \mathrm{sgn}(x_i)\, ,
\end{equation}
which depends only on the local state $x_i$ of agent $i$. We will exploit this fact in our reformulation of Problem \ref{prob:distributed_stress_minimization_approx}. 
The following lemma \textemdash{} whose proof may be found in Appendix~\ref{app:proof_lemma_strong_convexity} \textemdash{} characterizes the convexity of $f_{\alpha,\epsilon}(x)$.


\begin{lemma}[Strong convexity of $f_{\alpha,\epsilon}$]\label{lemma:strong_convexity}
Consider the function $f_{\alpha,\epsilon}:\R^{n_\ell}\mapsto\R$ defined as in \eqref{eq:cost_function}. Then, for all $0<\epsilon\leq 1$ and for all $\alpha>0$, $f_{\alpha,\epsilon}$ is strongly convex in $x$.
\end{lemma}


While Lemma \ref{lemma:recover_infty_norm} provides an asymptotic result, one finds that numerical issues are encountered for sufficient large (resp., small) values of $\alpha$ (resp., $\epsilon$). In practice however, reasonably small values of $\alpha$ (resp., $\epsilon$) provide excellent results, with no noticeable numerical issues.
\subsection{A Distributed Primal-Dual Feedback Controller}\label{subsec:primal_dual_alg} 
We now reformulate Problem \ref{prob:distributed_stress_minimization_approx} by exploiting \eqref{eq:cost_function} and the equivalence in \eqref{eq:argmin_equiv}. By introducing the additional quantities
\begin{equation}\label{eq:standard_constraints}
\mathcal{I} := \begin{bmatrix} \mathds{I} \\ -\mathds{I} \end{bmatrix}\, ,\ \ \ 
\chi := \begin{bmatrix} x_{\max} \\ -x_{\min}\end{bmatrix}\, ,\ \ \
\varphi := \begin{bmatrix} q_{\max} \\ -q_{\min}\end{bmatrix}\, ,\ \ \
\end{equation}
the smooth approximation of Problem \ref{prob:distributed_stress_minimization_approx} reads as
\begin{problem}[Smooth Stress Minimization]\label{prob:distributed_smooth_stress_minimization}
\begin{eqnarray}\label{eq:distributed_smooth_stress_minimization}
&	\underset{x\in\mathbb{R}^{n_\ell}}{\mathrm{minimize}}	& 	f_{\alpha,\epsilon}(x),\\
&	\mathrm{s.t.}		& 	\begin{cases}
							\mathcal{I} x \leq \chi\, ,\\
							\mathcal{I} (\Qcrit x - Q_L) \leq \varphi\, .
							\end{cases}							
							\nonumber							
\end{eqnarray}
\end{problem}
where, thanks to \eqref{eq:standard_constraints}, the constraints are now in the standard form of a convex program.\\
One possible way to solve Problem \ref{prob:distributed_smooth_stress_minimization} is the use of standard dual-ascent discrete-time algorithm \cite{DPB:08}. The Lagrangian function associated to \eqref{prob:distributed_smooth_stress_minimization} is equal to
%
\begin{equation*}
\resizebox{\hsize}{!}{
$\mathcal{L}(x,\lambda,\mu) = f_{\alpha,\epsilon}(x) + \lambda^T\left(\mathcal{I} x - \chi \right) + \mu^T\left(\mathcal{I}(\Qcrit x - Q_L) - \varphi \right)\,,$}
\end{equation*}
where $\lambda, \mu \in \mathbb{R}^{2n_{\ell}}$ are vectors of Lagrange multipliers. The dual-ascent then consists of the iterative updates
\begin{subequations}\label{eq:dual_ascent_dyn}
\begin{align}
x(t+1) &= \underset{x}{\mathrm{argmin}}\ \mathcal{L}(x,\lambda(t),\mu(t))\, ,\label{eq:dual_ascent_primal}\\
\lambda(t+1) &= \left[ \lambda(t) + \rho(\mathcal{I} x(t+1) - \chi)\right]^+\, ,\label{eq:dual_ascent_dual_lambda}\\
 \mu(t+1) &= \Big[ \mu(t) +\rho(\mathcal{I}(\Qcrit x(t+1) - Q_L) - \varphi)\Big]^+\label{eq:dual_ascent_dual_mu}
\end{align}
\end{subequations}
where $\rho > 0$ is the step size and $[\cdot]^+$ denotes the projection on the positive orthant, that is, $[a]^+ = a,$ if $a>0$, and $0$ otherwise.
%
%
We now state a convergence result whose proof can be found in Appendix \ref{app:proof_prop_dual_ascent_convergence}.

\smallskip 

\begin{proposition}[Convergence of the dual ascent algorithm]\label{prop:dual_ascent_convergence}
Consider Problem \ref{prob:distributed_smooth_stress_minimization} and assume that Slater's condition holds, namely, a strictly feasible solution for \eqref{eq:distributed_smooth_stress_minimization} exists. Then, there exists $\overline{\rho}$ such that for any $\rho\leq\overline{\rho}$, the dual ascent algorithm \eqref{eq:dual_ascent_dyn} converges to the optimal solution of \eqref{eq:distributed_smooth_stress_minimization}.
\end{proposition}

\smallskip

Observe that the desired control injections can be computed by the corresponding agents from \eqref{eq:control_input_distributed} as
\begin{equation}\label{eq:control_inputs}
q(t+1) = -(\Qcrit x(t+1) - Q_L)\, .
\end{equation}
From \eqref{eq:dual_ascent_dyn}, one can see that the proposed dual ascent algorithm is amenable of distributed implementation meaning that node $i$, to compute $x_i(t+1),\ \mu_i(t+1)$, needs only information coming from neighboring nodes (with the sparsity pattern induced by the $\Qcrit$ matrix); while no exchange of information is needed to compute $\lambda_i(t+1)$. To be more precise, observe that \eqref{eq:dual_ascent_primal} is separable in the $x_i$ variable; indeed, from the first order optimality condition must hold for any $i\in\mathcal{V}$ that
$$
\frac{\partial}{\partial x_i} \mathcal{L}(x,\lambda(t),\mu(t)) = 0\, ,
$$
which is equivalent to
\begin{align}
\frac{\partial f_{\alpha,\epsilon}(x)}{\partial x_i} = \underset{\zeta_i(t)}{\underbrace{-\sum_{j\in\mathcal{N}_i}\left([\mathcal{I}^T]_{ij}\lambda_j(t) + [\Qcrit\mathcal{I}^T]_{ij}\mu_j(t)\right)}}\, ,\label{eq:first_order_optimality}
\end{align}
%
%
\begin{figure}[t]
\centering
\includegraphics[width = \columnwidth]{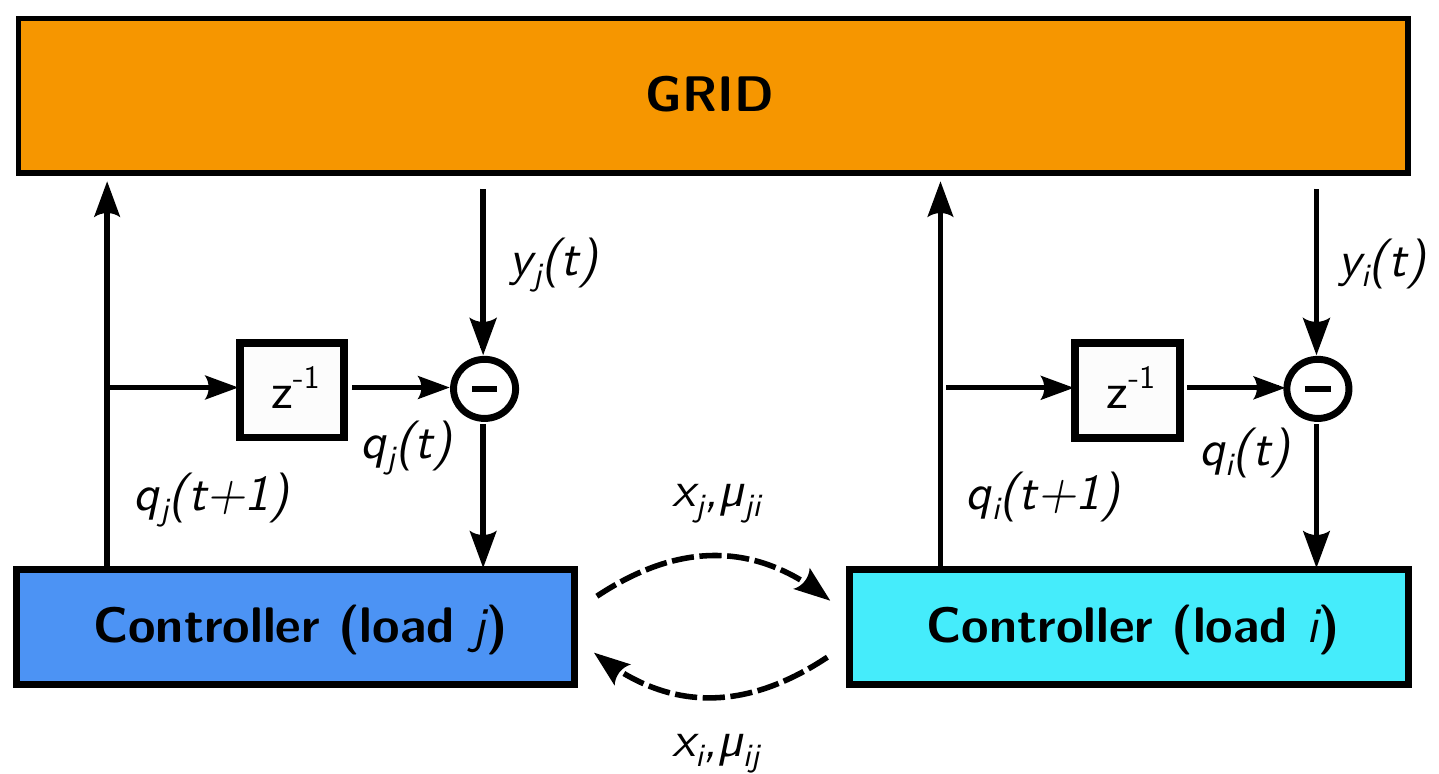}
\caption{\small 
Illustration of distributed control architecture.
Measurements $y(t)$ are gathered from the grid and used, together with $q(t)$, to update $x$. Neighboring controllers exchange  information. The block $z^{-1}$ is the one-step delay operator.}
\label{fig:control_loop}
\end{figure}
where $\mathcal{N}_i:=\{j\in\mathcal{V}\ : [\Qcrit]_{ij} \neq 0\}$ represents the set of neighbors of node $i$. Observe that the right hand side of \eqref{eq:first_order_optimality} is constant given the multipliers of the neighbors and that, thanks to strong convexity of $f_{\alpha,\epsilon}$ and  monotonicity of its first derivative, \eqref{eq:first_order_optimality} has always a unique real-valued solution. Interestingly, for $\epsilon=1$ a closed-form solution for \eqref{eq:first_order_optimality} exists and is equal to
\begin{equation}\label{eq:closed_form_primal}
x_i(t+1) = \mathrm{sgn}(\zeta_i(t)) \frac{1}{\sqrt{2\alpha}} \sqrt{W\left(\frac{\zeta_i^2(t)}{2\alpha}\right)}\, , 
\end{equation} 
where $W(\cdot)$ is the \emph{Lambert W} or \emph{ProductLog} function \cite{RMC-GHG-DEGH-DJJ-DEK:96} defined as the inverse of the function $g(z)=z\exp(z)$. Finally, note that, each agent needs to know some model information, namely its corresponding $\Qcrit$ entries which are related to the electric quantities connecting them to their neighbors.\\
It is worth noticing the presence of the load demand $Q_L$ in both \eqref{eq:dual_ascent_dual_mu} and \eqref{eq:control_inputs}. Usually, it is reasonable to have voltage and current monitoring at each bus. From this measurements it is possible to extract information about the total reactive load absorbed or injected at the bus. In particular, by defining $y:= Q_L + q$ as the aggregate contribution of the load together with the control input, from the bus monitoring it is possible to measure $y$ rather than $Q_L$ independently. In order to compute \eqref{eq:dual_ascent_dual_mu} and \eqref{eq:control_inputs} it is sufficient to set
$$
Q_L = y(t) - q(t)\, ,
$$
where $y(t)$ and $q(t)$ are the aggregate load measurements and the control input at the current $t$ iteration. An illustration of the control loop is showed in Figure \ref{fig:control_loop}: the aggregate measurements are taken from the grid and are used together with the control input to compute the update $x(t+1)$. Then, by using $x(t+1)$, $y(t)$ and $q(t)$, the new control input $q(t+1)$, used to actuate the grid, are computed.
%
\begin{figure}[t]
\centering
\includegraphics[width=\columnwidth]{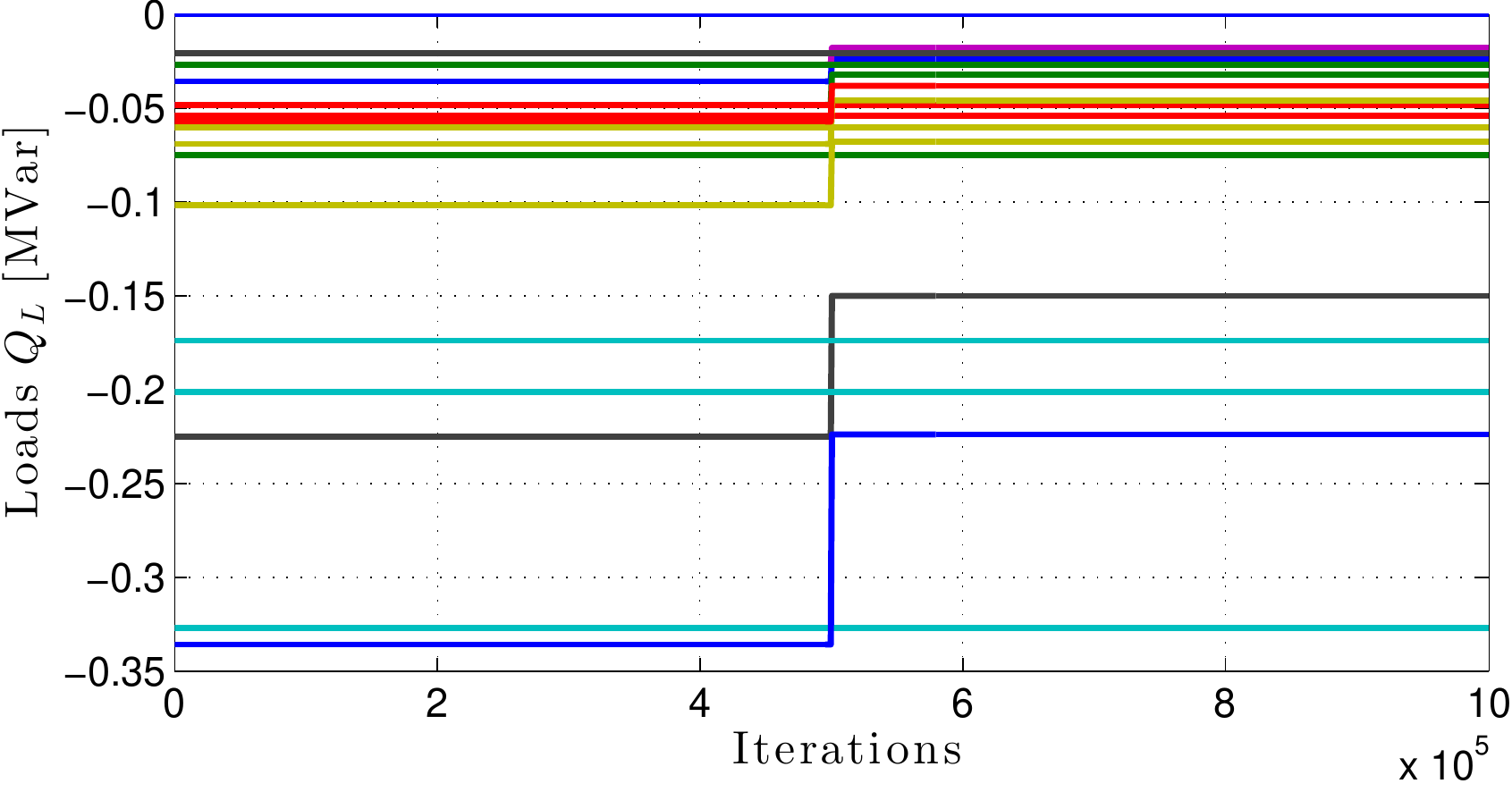}
\caption{\small Dynamics of the loads $Q_L$. The different colored curves represent the time evolution of the loads at different buses.}
\label{fig:load_dyn}
\end{figure}
\begin{figure}[b]
\centering
\includegraphics[width = \columnwidth]{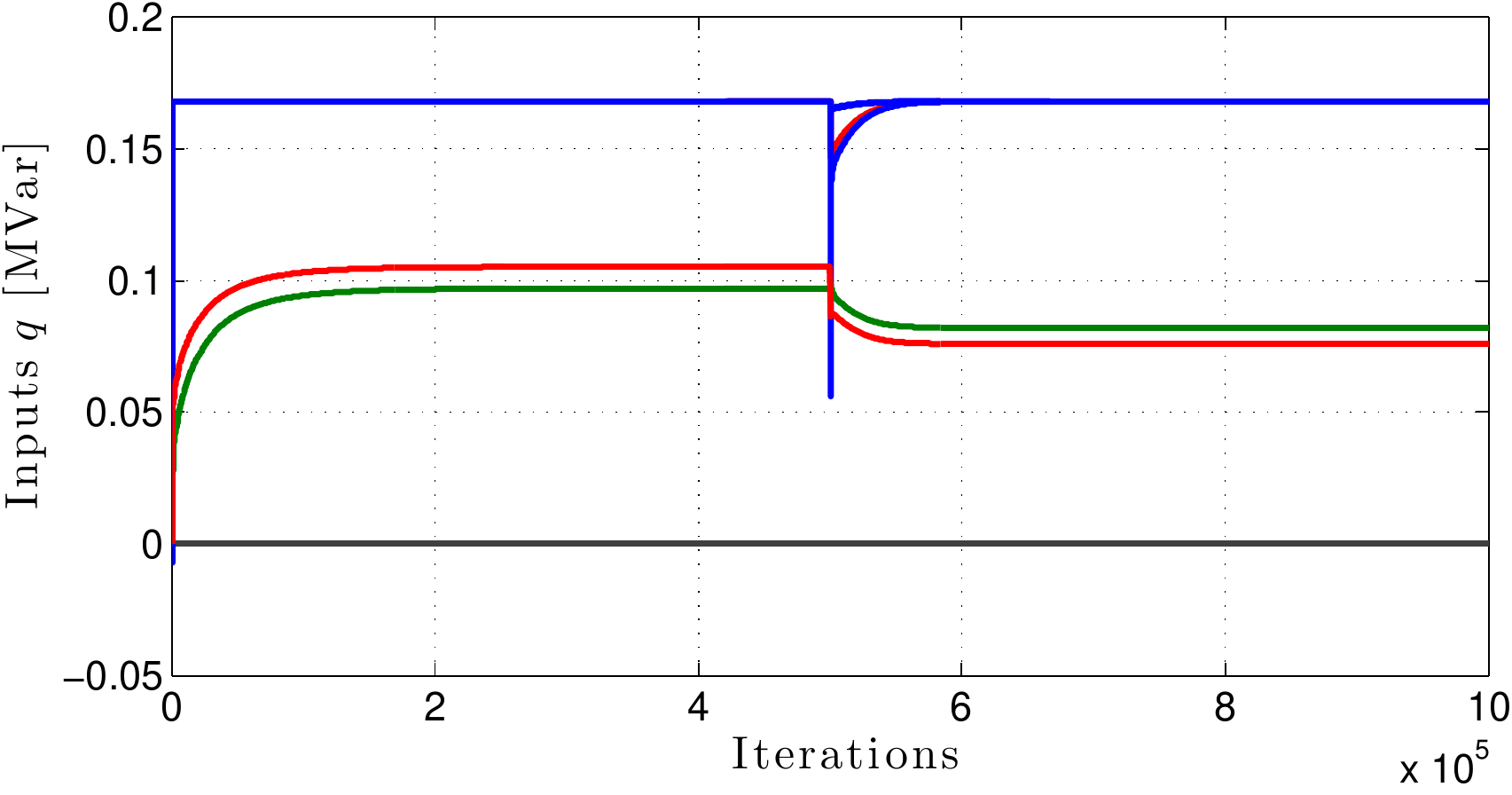}
\caption{\small Dynamic of the control inputs $q$. The different colored curves represent control inputs at different controlled buses.}
\label{fig:injection_dyn}
\end{figure}
\subsection{Simulation: Distributed Online Feedback Controller}\label{subsec:simulation_distributed} 
\begin{figure}[t]
\centering
\includegraphics[width = \columnwidth]{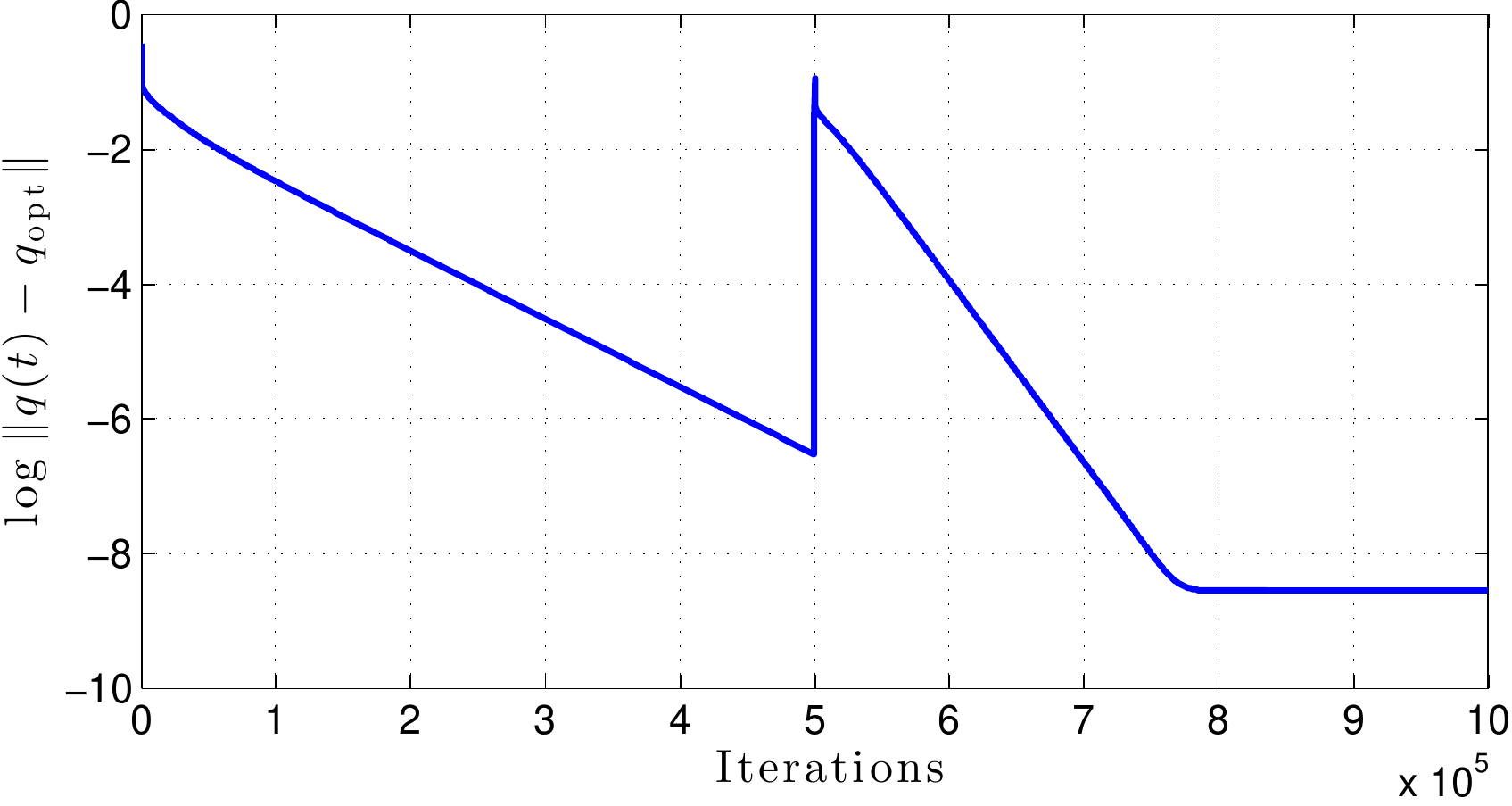}
\caption{\small Evolution of the norm of the error between the online values of $q(t)$ and their optimal value $q_{\rm opt}$ computed offline.}
\label{fig:distanceToOptimalInjection}
\end{figure}

We now present simulations illustrating the effectiveness of our distributed controller in the presence of time-varying loads. The setup is the same as used in Section \ref{subsec:simulation_sparsity}. However, differently from the previous set of simulations, we simulate the full coupled system at each iteration of the algorithm by means of MATPOWER \cite{RDZ-CEMS-RJT:11}. We show the performance of the algorithm assuming we have access to aggregate measurements $y$ for a fixed allocation of the resources. Specifically, there are 6 controllable units out of 24 total loads at the load buses number 1, 10, 15, 22, 23 and 24. Moreover, the controlled injections $q$ are saturated before actuating the electric grid, simulating the fact that they cannot exceed their capacity limits. The loads (active and reactive) randomly encounter, at half of the simulation time, a jump equal to $40\%$ of their starting value (see Figure~\ref{fig:load_dyn}). All other parameters are held constant. Figure~\ref{fig:injection_dyn} shows the evolution of the reactive injection returned by the algorithm.
%
\begin{figure}[t]
\centering
\includegraphics[width = \columnwidth]{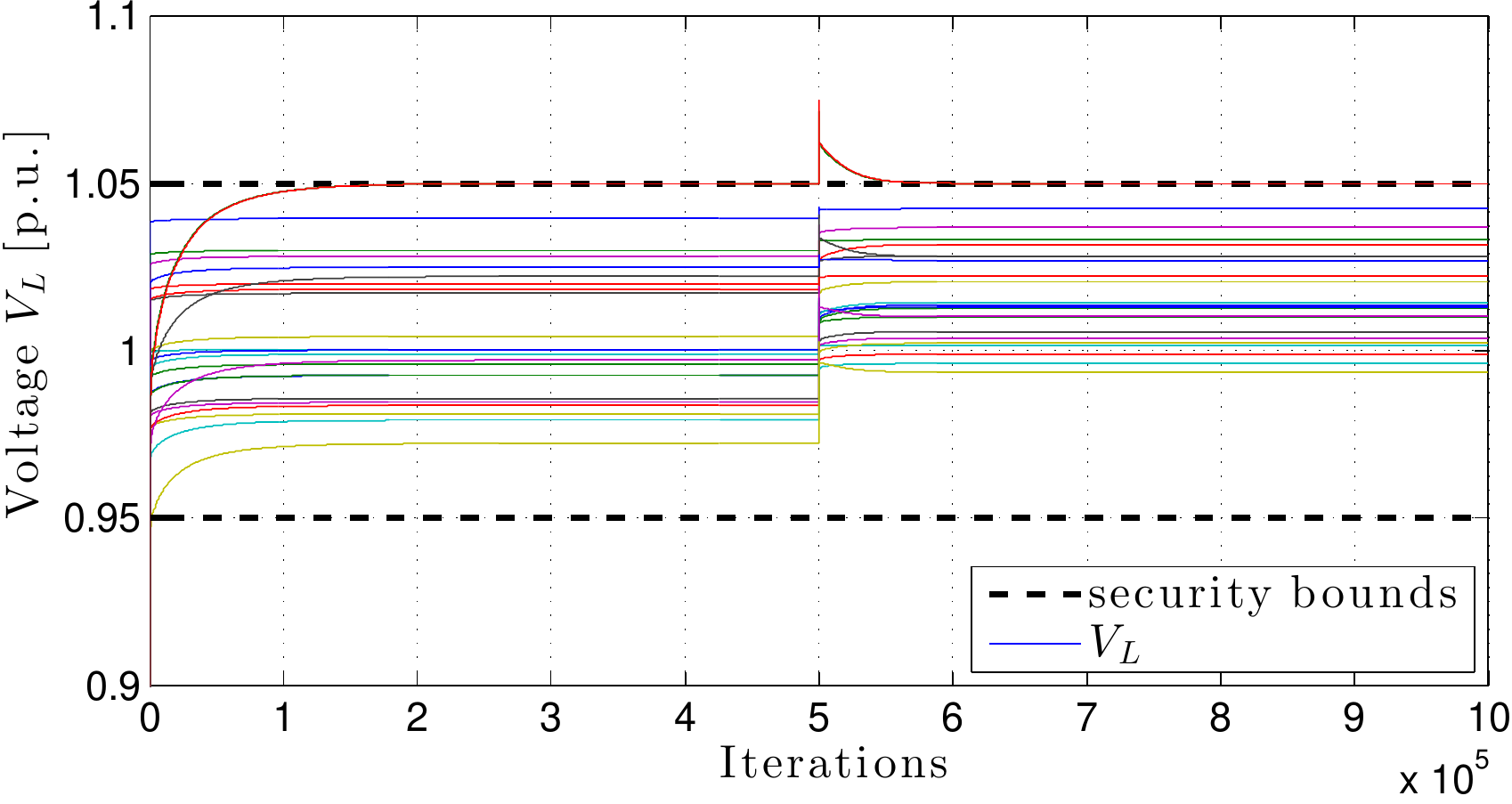}
\caption{\small Evolution of the voltages corresponding to the coupled nonlinear power flow equations.}
\label{fig:voltages_nonlinearPF}
\end{figure}
%
Figure \ref{fig:voltages_nonlinearPF} shows the corresponding evolution of bus voltages under the distributed controller \eqref{eq:dual_ascent_dyn}--\eqref{eq:control_inputs}. Observe that the voltage solution of the full coupled PFEs remains within the operational bounds. Moreover, we note that the worst case steady-state difference between the solution of the full coupled PFEs and the linearized solution of the RPFEs given by \eqref{eq:LoadRPFEnormalize_controlled} is of only $1.2\%$. This fact highlights the effectiveness of both the linearization and the control algorithm. Finally, Figure \ref{fig:distanceToOptimalInjection} shows the evolution of the error between the values of the injection $q(t)$ and the optimal value computed offline $q_{\rm opt}$ as solution of Problem \ref{prob:stress_minimization_approx}, as a function of the iterations, in logarithmic scale. Notice that the value of $q$ computed online converges to the optimal offline value even after the change in the loads. As a limit, it must be noted that the proposed distributed algorithm requires substantial amount of iteration until convergence.
\section{Conclusions and Future Directions}\label{sec:conclusions}
We considered the problem of voltage support via reactive power injections. Conversely to what suggested by conventional wisdom, we showed that the standard local security requirement might be inadequate. Then, we proposed a novel optimization formulation whose cost function encodes the stress experienced by the grid, while the security requirements are imposed as constraints. Thanks to recent advances on the solvability and linearization of the reactive power flows, the problem becomes linear and convex and can be efficiently solved. In addition, we addressed the planning problem to solve for the optimal allocation of the resources. Finally, thanks to a suitable reformulation, we presented a distributed algorithm to solve for the stress minimization, which implements a real-time feedback controller and which, as drawback, is characterized by a slow convergence rate.

As future research, relevant practical applications would be to seek for distributed implementation which are faster in the convergence rate and employ only communications among compensators. Moreover, it would be worth exploring our novel formulation, which combines system-level stress minimization and local security constraints, with respect to different optimization variables. From a power system analysis perspective the ultimate problem remains the analysis of the full coupled power flow equations.
%
%

\appendix

\subsection{Proof of Lemma \ref{lemma:recover_infty_norm}}\label{app:proof_lemma_recover_infty_norm}
In the following, we will prove the result by first taking the limit for $\alpha$ followed by the limit for $\epsilon$. It is easy to show, thanks to a Taylor series expansion of $\tilde{f}_{\alpha,\epsilon}$ around $\epsilon=0$, that exchanging the order of the limits does not change the result.\\ 
Consider the function $\tilde{f}_{\alpha,\epsilon}(x)$ as defined in \eqref{eq:smooth_infty} and 
let us define $|x|_{\max}^{1+\epsilon}:= \max_i |x_i|^{1+\epsilon}$. It is possible to rewrite
$$
\tilde{f}_{\alpha,\epsilon}(x) = |x|_{\max}^{1+\epsilon} + \frac{1}{\alpha}\log\Big(\frac{1}{n}\sum_i\exp\Big(\alpha(\underset{\leq 0}{\underbrace{|x_i|^{1+\epsilon} - |x|_{\max}^{1+\epsilon}}})\Big)\Big)\, .
$$
Now, since the exponent in the second term is always $<0$ except for the components where the maximum is attained which are exactly equal to $0$, we have
$$
\lim_{\alpha\rightarrow+\infty} \tilde{f}_{\alpha,\epsilon} = |x|_{\max}^{1+\epsilon}\, .
$$
Finally, to conclude the proof, notice that the Taylor series expansion around $\epsilon=0^+$ of $|x|^{1+\epsilon}_{\max}$ is equal to
$$
|x|_{\max} + \sum_{n=1}^{\infty}\frac{1}{n!} |x|_{\max}\log^n(|x|_{\max}) \epsilon^n \underset{\epsilon\rightarrow0^+}{\longrightarrow} |x|_{\max}\, .\quad\quad\oprocend
$$

\subsection{Proof of Lemma \ref{lemma:strong_convexity}}\label{app:proof_lemma_strong_convexity}
To prove strong convexity of $f_{\alpha,\epsilon}(x)$ we exploit the second order characterization of strong convexity which states that a function is strongly convex if and only if 
$
\nabla^2_xf_{\alpha,\epsilon}(x) - m\mathds{I}\ 
$
is positive definite, 
for some $m>0$. The Hessian of $f_{\alpha,\epsilon}$ is indeed a diagonal matrix whose $i$-th diagonal entry is equal to
$$
\underset{> 0}{\underbrace{\alpha(1+\epsilon)\exp(\alpha|x_i|^{1+\epsilon})}}\left( \epsilon \frac{|x_i|^{\epsilon}}{|x_i|} + \alpha(1+\epsilon)|x_i|^{2\epsilon}\right).
$$
Note that, the entire expression can fail to be positive only if the second term in the right hand side possibly fails to be positive. The only possible point of failure is $x_i=0$ where we encounter a $\frac{0}{0}$ limit. However, being $|x_i|^\epsilon = o(|x_i|)$ for any finite value $0<\epsilon\leq 1$, we have that
$$
\frac{\partial^2 f_{\alpha,\epsilon}(x)}{\partial x_i^2} \underset{|x_i|\rightarrow 0}{\longrightarrow} +\infty\, .
$$
Since each diagonal component is bounded below by $m_i>0$, defining $m:=\min_i\ m_i$, we can conclude.
\oprocend

\subsection{Proof of Proposition \ref{prop:dual_ascent_convergence}}\label{app:proof_prop_dual_ascent_convergence}
First of all, we recall (assuming Slater's condition holds) that the duality gap is zero. Hence, solving the primal Problem \ref{prob:distributed_smooth_stress_minimization} is equivalent to solve its dual problem which is defined as
\begin{eqnarray}
&	\underset{\lambda}{\mathrm{maximize}}	& 	d(\lambda) := \inf_x\ f_{\alpha,\epsilon}(x) + \lambda^T\Big(A x - b\Big)\, ,\label{eq:dual_problem}\\
&	\mathrm{s.t.}		& 	\lambda \geq 0\, .	\notag						
\end{eqnarray}
where $A:= \begin{bmatrix}\mathcal{I}^T & (\mathcal{I}\Qcrit)^T\end{bmatrix}^T$ and $b:=\begin{bmatrix}\chi^T & \varphi^T\end{bmatrix}^T$.
Moreover, thanks to strong convexity of $f_{\alpha,\epsilon}$, the solution is unique. From Proposition $6.1.1$ in \cite{DPB:08} follows that $d$ is everywhere continuously differentiable. Moreover
$$
\nabla_\lambda d(\lambda) = Ax^*_\lambda - b
$$
where $x^*_\lambda := \mathrm{argmin}_x\ f_{\alpha,\epsilon}(x) + \lambda^Tg(x)$. In addition, the dual ascent algorithm \eqref{eq:dual_ascent_dyn} coincides with a projected gradient \cite{DPB:08} applied to \eqref{eq:dual_problem} which, thanks to Proposition $2.3.2$ in \cite{DPB:08}, is known to converge, for sufficiently small step sizes, namely $0 < \rho < \frac{2}{L} = \overline{\rho}$, if $\nabla_\lambda d$ is Lipschitz continuous with Lipschitz constant $L$. To prove Lipschitz continuity of $\nabla_\lambda d$, observe that it is linear respect to $x^*_\lambda$. Then, we just need Lipschitz continuity of $x_\lambda^*$ respect to $\lambda$. From the definition of $x^*_\lambda$ and thanks to first order optimality condition, it holds that
$$
\nabla_xf_{\alpha,\epsilon}(x_\lambda^*) = -A^T\lambda\, .
$$
By defining $\zeta:=-A^T\lambda$ and recalling from \eqref{eq:cost_derivative} that each component of $\nabla_x f$ is only a function of $x_i$, we have that 
$$
\frac{\partial f_{\alpha,\epsilon}(x^*_\lambda)}{\partial x_i} = \zeta_i\, ,
$$
and it is possible to reduce the analysis to show Lipschitz continuity of the $i$-th component of $x^*_\lambda$ respect to $\zeta_i$. Now, by being $f_{\alpha,\epsilon}$ twice continuously differentiable and thanks to the inverse function theorem, $\partial f_{\alpha,\epsilon}/\partial x_i$ is invertible, namely 
$$
[x_\lambda^*]_i = \left(\frac{\partial f_{\alpha,\epsilon}}{\partial x_i}\right)^{-1}(\zeta_i)\, ,
$$
its inverse is continuously differentiable, and moreover
$$
\frac{\partial \left((\partial f_{\alpha,\epsilon}/\partial x_i)^{-1}\right)}{\partial x_i} =  \left(\frac{\partial(\partial f_{\alpha,\epsilon}/\partial x_i)}{\partial x_i}\right)^{-1} < \frac{1}{m_i}\, ,
$$
where the last inequality holds since $f_{\alpha,\epsilon}$ is strongly convex. Then, $x^*_\lambda$ is Lipschitz continuous respect to $\zeta$ and so respect to $\lambda$ being the former a linear function of the latter.
\oprocend

\bibliographystyle{IEEEtran}



\bibliography{/Users/marcotodescato/Documents/Work/repository/Svn/bib/alias,/Users/marcotodescato/Documents/Work/repository/Svn/bib/Main,/Users/marcotodescato/Documents/Work/repository/Svn/bib/New,/Users/marcotodescato/Documents/Work/repository/Svn/bib/FB}



\end{document}